\documentclass[leqno,a4paper]{article}

\usepackage{a4wide}

\usepackage{amssymb,amsfonts,amsthm,mathrsfs}

\usepackage[sumlimits,intlimits,namelimits]{amsmath}

\usepackage{epsfig}

\usepackage{pst-all}
\usepackage{multido}

\psset{xunit=1cm,yunit=1cm,runit=1cm}

\newcommand{\csub}{
 \subset\!\subset
}

\newcommand{\Cinf}{\ensuremath{\mathcal{C}^\infty}}

\newcommand{\D}{\ensuremath{{\cal D}}}

\newcommand{\loc}{\ensuremath{\text{loc}}}

\newcommand{\mb}[1]{\ensuremath{\mathbb{#1}}}
\newcommand{\N}{\mb{N}}

\newcommand{\R}{\mb{R}}
\newcommand{\C}{\mb{C}}

\newcommand{\cl}[1]{\ensuremath{[#1]}}

\newcommand{\G}{\ensuremath{{\cal G}}}

\newcommand{\Gc}{\ensuremath{{\cal G}_\mathrm{c}}}
\newcommand{\EM}{\ensuremath{{\cal E}_{\mathrm{M}}}}

\newcommand{\NN}{\ensuremath{{\cal N}}}

\renewcommand{\d}{\ensuremath{\partial}}
\newcommand{\diff}[1]{\frac{d}{d#1}}

\renewcommand{\div}{\ensuremath{\mbox{\rm div}\,}}

\newfont{\bl}{msbm10 scaled \magstep2}

\newtheorem{theorem}{Theorem}[section]
\newtheorem{lemma}[theorem]{Lemma}

\newtheorem{definition}[theorem]{Definition}

\theoremstyle{definition}
\newtheorem{remark}[theorem]{Remark}
\newtheorem{example}[theorem]{Example}

\newcommand{\beq}{\begin{equation}}
\newcommand{\eeq}{\end{equation}}

\newcommand{\emb}{\hookrightarrow}

\newcommand{\col}{\colon}

\newcommand{\dis}[2]{\langle #1 , #2 \rangle}
\newcommand{\inp}[2]{\langle #1 | #2 \rangle}  
\newcommand{\notmid}{\mid\kern-0.5em\not\kern0.5em}

\newcommand{\norm}[2]{{\| #1 \|}_{#2}}

\newcommand{\Norm}[1]{\norm{#1}{}}

\newcommand{\al}{\alpha}

\newcommand{\eps}{\varepsilon}

\newcommand{\vphi}{\varphi}

\newcommand{\Om}{\Omega}

\renewcommand{\Re}{\ensuremath{\mathop{\mathrm{Re}}}}

\newcommand{\ovl}[1]{\overline{#1}}

\begin{document}

\title{{\bf Comparison of some solution concepts for linear first-order hyperbolic differential equations with non-smooth coefficients}\footnote{Work supported by FWF grant
Y237-N13}}

\author{
\emph{Simon Haller}\footnote{\texttt{simon.haller@univie.ac.at}}  \enspace \& \enspace 
\emph{G\"{u}nther H\"{o}rmann}\footnote{\texttt{guenther.hoermann@univie.ac.at}} \\
\ \\
\ \\ 
Fakult\"at f\"ur Mathematik\\
Universit\"at Wien, Austria \\
}

\date{\today}
\maketitle

\vspace{-1cm}
\begin{abstract}
We discuss solution concepts for linear hyperbolic equations with coefficients of regularity below Lipschitz continuity. Thereby our focus is on theories which are based either on a generalization of the method of characteristics or on refined techniques concerning energy estimates. We provide a series of  examples both as simple illustrations of the notions and conditions involved but also to show logical independence among the concepts.\\[1mm] 
\end{abstract}

\setcounter{section}{-1}

\section{Introduction}

According to Hurd and Sattinger in \cite{HS:68}
the issue of a systematic investigation of hyperbolic partial differential equations with discontinuous coefficients as a research topic has been raised by Gelfand in 1959. Here, we attempt a comparative study of some of the theories on that subject which have been put forward since. More precisely, we focus on techniques and concepts that build either on the geometric picture of propagation along \emph{characteristics} or on the functional analytic aspects of \emph{energy estimates}. 

In order to produce a set-up which makes the various methods comparable at all, we had to stay with the special situation of a \emph{scalar partial differential equation with real coefficients}. As a consequence, for example, we do not give full justice to theories whose strengths lie in the application to systems rather than to a single equation. 
A further limitation in our choices comes from the restriction to concepts, hypotheses and mathematical structures which (we were able to) directly relate to distribution theoretic or measure theoretic notions.

To illustrate the basic problem in a simplified lower dimensional situation for a linear conservation law, we consider the following formal differential equation for a density function (or distribution, or generalized function) $u$ depending on time $t$ and spatial position $x$ 
$$
   \d_t u(t,x) + \d_x(a(t,x) u(t,x)) = 0.
$$
Here, $a$ is supposed to be a \emph{real} function (or distribution, or generalized function) and the derivatives shall be interpreted in the distributional or weak sense. This requires either to clarify the meaning of the product $a \cdot u$ or to avoid the strict meaning of ``being a solution''.  

An enormous progress has been made in research on nonlinear conservation laws (cf., e.g.\ \cite{Hoermander:97,BGS:07} and references therein)  of the form
$$
  \d_t u(t,x) + \d_x( g(u(t,x)) ) = 0,  
$$
where $g$ is a (sufficiently) smooth function and $u$ is such that $g(u)$ can be defined in a suitable Banach space of distributions. Note however, that this equation  does not include linear operators of the form described above as long as the nonlinearity $g$ does not include
additional dependence on $(t,x)$ as independent variables (i.e., is not of the more general form $g(t,x,u(t,x))$). Therefore the theories for linear equations described in the present paper are typically not mere corollaries of the nonlinear theories.
Essentially for the same reason we have also not included methods based on Young measures (cf.\ \cite[Chapter V]{Hoermander:97}).

Further omissions in our current paper concern hyperbolic equations of second order. For advanced theories on these we refer to the energy method developed by Colombini-Lerner in \cite{CL:95}. An overview and illustration of non-solvability or non-uniqueness effects with wave equations and remedies using Gevrey classes can be found in \cite{Spagnolo:87}.

Of course, also the case of first-order equations formally ``of principal type'' with non-smooth complex coefficients is of great interest. It seems that the borderline between solvability and non-solvability is essentially around Lipschitz continuity of the coefficients (cf.\ \cite{Jacobowitz:92,Hounie:91,HounieMelo:95}). Moreover, the question of uniqueness of solutions in the first-order case has been addressed at impressive depth in \cite{CL:02}.

Our descriptive tour with examples consists of two parts: Section 1 describes concepts and theories extending the classical method of characteristics, while Section 2 is devoted to theories built on energy estimates. All but two of the theories or results (namely, in Subsections 1.3 and 2.3.2) we discuss and summarize are not ours. However, we have put some effort into unifying the language and the set-up, took care to find as simple as possible examples which are still capable of distinguishing certain features, and have occasionally streamlined or refined the original or well-known paths in certain details. 

In more detail, Subsection 1.1 starts with Caratheodory's theory of generalized solutions to first-order systems of (nonlinear) ordinary differential equations and adds a more distribution theoretic view to it. In Subsection 1.2 we present the generalization in terms of Filippov flows and the application to transport equations according to Poupaud-Rascle. Subsection 1.3 provides a further generalization of the characteristic flow as Colombeau generalized map with nice compatibility properties when compared to the Filippov flow. In Subsection 1.4 we highlight some aspects or examples of semigroups of operators on Banach spaces stemming from underlying generalized characteristic flows on the space-time domain. We also describe a slightly exotic concept involving the measure theoretic adjustment of coefficients to prescribed characteristics  for $(1+1)$-dimensional equations according to Bouchut-James in Subsection 1.5.

Subsection 2.1 presents a derivation of energy estimates under very low regularity assumptions on the coefficients and also discusses at some length the functional analytic machinery to produce a solution and a related weak solution concept for the Cauchy problem. Subsection 2.2 then compares those three theories, namely by Hurd-Sattinger, Di Perna-Lions, and Lafon-Oberguggenberger,  which are based on regularization techniques combined with energy estimates. Finally, Subsection 2.3 briefly describes two related results obtained by paradifferential calculus, the first concerning energy estimates and the solution of the Cauchy problem for a restricted class of operators, the second is a method to reduce equations to equivalent ones with improved regularity of the source term.

As it turns out in summary, none of the solution concepts for the hyperbolic partial differential equation is contained in any of the others in a strict logical sense. However, there is one feature of the Colombeau theoretic approach: it is always possible to model the coefficients and initial data considered in any of the other theories (by suitable convolution regularization) in such a way that the corresponding Cauchy problem becomes uniquely solvable in Colombeau's generalized function algebra.  In many cases the Colombeau generalized solution can be shown to have the appropriate distributional aspect in the sense of heuristically reasonable solution candidates.


\subsection{Basic notation and spaces of functions, distributions, and generalized functions}

Let $\Omega$ denote an open subset of $\R^n$. We use the notation $K \Subset \Omega$, if $K$ is a compact subset of $\Omega$. The letter $T$ will always be used for real number such that $T > 0$. We often write $\Omega_T$ to mean $]0,T[ \times \R^n$ with closure $\ovl{\Omega_T} = [0,T] \times \R^n$.

The space $C^\infty(\ovl{\Omega})$ consists of smooth functions on $\Omega$ all whose derivatives have continuous extensions to $\ovl{\Omega}$. For any $s \in \R$ and $1 \leq p \leq \infty$ we have the Sobolev space $W^{s,p}(\R^n)$ (such that $W^{0,p} = L^p$), in particular $H^s(\R^n) = W^{s,2}(\R^n)$. Our notation for $H^s$-norms and inner products will be ${\Norm{.}}_s$ and $\dis{.}{.}_s$, in particular, this reads ${\Norm{.}}_0$ and $\dis{.}{.}_0$ for the standard $L^2$ notions.

We will also make use of the variants of Sobolev and $L^p$ spaces of functions on an interval $J \subseteq \R$ with values in a Banach space $E$, for which we will employ a notation as in $L^1(J;E)$, for example. (For a compact treatment of the basic constructions we refer to \cite[Sections 24 and 39]{Treves:75}.) Furthermore, as usually the subscript 'loc' with such spaces will mean that upon multiplication by a smooth cutoff we have elements in the standard space. We occasionally write $AC(J;E)$ instead of $W^{1,1}_{\rm{loc}}(J;E)$ to emphasize the property of absolute continuity. 

The subspace of Distributions of order $k$ on $\Omega$ ($k \in \N$, $k \geq 0$) will be denoted by $\D'^k(\Omega)$. We identify $\D'^0(\Omega)$ with the space of complex Radon measures $\mu$ on $\Omega$, i.e., $\mu = \nu_+ - \nu_- + i(\eta_+ - \eta_-)$, where $\nu_\pm$ and $\eta_\pm$ are positive Radon measures on $\Omega$, i.e., locally finite (regular) Borel measures. 

As an alternative regularity scale with real parameter $s$ we will often refer to the H\"older-Zygmund classes $C_\ast^s(\R^n)$  (cf.\ \cite[Section 8.6]{Hoermander:97}). In case $0 < s < 1$ the corresponding space comprises the continuous bounded functions $u$ such that there is $C > 0$ with the property that for all $x \neq y$ in $\R^n$ we have
$$
   \frac{|u(x) - u(y)|}{|x - y|^s} \leq C.
$$

Special types of distributions on $\R$ will be used in several of our examples to follow: the Heaviside function will be understood to be the $L^\infty(\R)$ class of the function defined almost everywhere by $H(x) = 0$ ($x < 0$), $H(x) = 1$ ($x > 0$), and will again be denoted by $H$; the signum function is $\rm{sign}(x) = H(x) - H(-x)$ ; furthermore, $x_+$ denotes the continuous function with values $x_+ = 0$ ($x < 0$), $x_+ = x$ ($x \geq 0$), $x_- = x_+ - x$; $\delta$ denotes the Dirac (point) measure at $0$ (in any dimension). 
\paragraph{Model product of distributions:}

A whole hierarchy of coherent distributional products has been discussed in
\cite[Chapter II]{O:92}, each of these products yielding the classical pointwise  multiplication when both factors are smooth functions. The most general level of this hierarchy is that of the so-called \emph{model product} of distributions $u$ and $v$, denoted by ${[u \cdot v]}$ if it exists.

We first regularize both factors by convolution with a \emph{model delta net} $(\rho_\eps)_{\eps > 0}$, where $\rho_\eps(x) = \rho(x/\eps)/\eps^n$ with $\rho\in\D(\R^n)$ such that $\int\rho(x)\, dx = 1$. Then the product of the corresponding smooth regularizations may or may not converge in $\D'$. If it does, the model product is defined by 
$$ 
  {[u \cdot v]} = \lim\limits_{\eps\to 0} (u*\rho_\eps) (v*\rho_\eps).
$$
In this case, it can be shown that the limit is independent of
the choice of $\rho$. For example, we have ${[H \cdot \delta]} = \delta /2$ and ${[\delta \cdot \delta]}$ does not exist.

\paragraph{Colombeau generalized functions:}

Our standard references for
the foundations and some applications of Colombeau's   nonlinear theory of generalized functions are \cite{Colombeau:84,Colombeau:85,O:92,GKOS:01}. 
We will employ the so-called special variant of Colombeau algebras, denoted by $\G^s$ in \cite{GKOS:01},
although here we shall simply use the letter $\G$ instead.

Let us briefly recall the basic constructions and properties. \emph{Colombeau generalized functions} on
$\Om$ are defined as equivalence classes $u = \cl{(u_\eps)_\eps}$ of nets of
smooth functions $u_\eps\in\Cinf(\Om)$ (\emph{regularizations}) subjected to
asymptotic norm conditions with respect to $\eps\in (0,1]$ for their
derivatives on compact sets: in more detail, we have
\begin{trivlist}
 
\item[$\bullet$] moderate nets $\EM(\Om)$: $(u_\eps)_\eps\in\Cinf(\Om)^{(0,1]}$ such that
 for all $K \Subset \Om$ and $\al\in\N^n$ there exists $p \in \R$ such that
 \beq\label{basic_estimate}
    \norm{\d^\al u_\eps}{L^\infty(K)} = O(\eps^{-p}) \qquad (\eps \to 0);
 \eeq

\item[$\bullet$] negligible nets $\NN(\Om)$: $(u_\eps)_\eps\in \EM(\Om)$ such that for all
  $K \Subset \Om$ and for all $q\in\R$ an
 estimate $\norm{u_\eps}{L^\infty(K)} = O(\eps^{q})$ ($\eps \to 0$) holds;

\item[$\bullet$] $\EM(\Om)$ is a differential algebra with operations defined at fixed $\eps$,
 $\NN(\Om)$ is an ideal, and $\G(\Om) := \EM(\Om) / \NN(\Om)$ is the (special)
 \emph{Colombeau algebra};
 
\item[$\bullet$] there are embeddings, $\Cinf(\Om)\emb \G(\Om)$ as a subalgebra and
 $\D'(\Om) \emb \G(\Om)$ as a linear subspace, commuting with partial derivatives;
 
\item[$\bullet$]  $\Om \to \G(\Om)$ is a fine sheaf and $\Gc(\Om)$ denotes
 the subalgebra of elements with compact support; by a
 cut-off in a neighborhood of the support one can always obtain representing nets with supports contained in a joint compact set;

\item[$\bullet$] in much the same way, one defines the Colombeau algebra $\G(\ovl{\Omega})$ on the closure of the open set $\Omega$ using representatives which are moderate nets in $C^\infty(\ovl{\Omega})$ (estimates being carried out on compact subsets of $\ovl{\Omega}$); 

\item[$\bullet$] two Colombeau functions $u = \cl{(u_\eps)_\eps}$ and $v = \cl{(v_\eps)_\eps}$ are said to be \emph{associated}, we write $u \approx v$, if $u_\eps - v_\eps \to 0$ in $\D'$ as $\eps \to 0$; furthermore, we call $u$ associated to the distribution $w \in \D'$, if $u_\eps \to w$ in $\D'$ as $\eps \to 0$; $w$ is then called the \emph{distributional shadow} of $u$ and we also write $u \approx w$;

\item[$\bullet$] assume that $\Omega$ is of the form $\Omega = \, ]T_1,T_2[ \, \times \Omega'$, where $\Omega' \subseteq \R^n$ open and $-\infty \leq T_1 < T_2 \leq \infty$; then we may define the restriction of $u = \cl{(u_\eps)_\eps} \in \G(\Omega)$ to the hyperplane $\{t_0\} \times \Omega'$ ($T_1 < t_0 < T_2$) to be the element $u \mid_{t = t_0} \G(\Omega')$ defined by the representative $(u_\eps(t_0,.))_\eps$; similarly, we may define the restriction of $u \in \G([T_1,T_2] \times \ovl{\Omega'})$ to $t = t_0$ for $T_1 \leq t_0 \leq T_2$ and obtain an element $u \mid_{t = t_0} \G(\ovl{\Omega'})$.

\item[$\bullet$] the set $\widetilde{\R}$ of \emph{Colombeau generalized real numbers} is defined in a similar way via equivalence classes $r = \cl{(r_\eps)_\eps}$ of nets of real numbers $r_\eps\in\R$ subjected to moderateness conditions $|r_\eps| = O(\eps^{-p})$ ($\eps \to 0$, for some $p$) modulo negligible nets satisfying $|r_\eps| = O(\eps^{q})$ ($\eps \to 0$, for all $q$); if $A \subset \R$ we denote by $\widetilde{A}$ the set of all generalized numbers having representatives contained in $A$ (for all $\eps \in ]0,1]$). Similarly, if $B \subset \R^n$ we construct $\widetilde{B} \subset {\widetilde{\R}}^n$ from classes of nets $(x_\eps)_\eps$ with $x_\eps \in B$ for all $\eps$;

\item[$\bullet$] a Colombeau generalized function $u = \cl{(u_\eps)_\eps} \in \G(\Omega)^d$ is said to be \emph{c-bounded} (compactly bounded), if for all $K_1 \Subset \Omega$ there is $K_2 \Subset \R^d$ and $\eps_0 > 0$ such that $u_\eps(K_1) \subseteq K_2$ holds for all $\eps > \eps_0$. 

\end{trivlist}

%
%
%
%
%
%
%
%
%
%
%
%
%
%
%

%
\section{Solution concepts based on the characteristic flow}
In this section we introduce solution concepts for first order partial differential equations,
which are based on solving the system of ordinary differential equations for the characteristics and using the resulting characteristic flow to define a solution.
To illustrate the basic notions we consider the following special case of the Cauchy problem in conservative form
$$
  L u :=\partial_t u + \sum_{k=1}^n \d_{x_k}( a_k(t,x)  u) = 0,
  \quad u(0) = u_0 \in\mathcal{D}'(\mathbb{R}^n),
$$ 
where the coefficients $a_k$ are real-valued bounded smooth functions. 
The associated system of ordinary differential equations for the characteristic curves reads
\begin{eqnarray*}
  \dot{\xi_k}(s) = a_k(s, \xi(s)) , \ \ \xi_k(t)=x_k \qquad 
  (k=1,\ldots,n).
\end{eqnarray*}
We use the notation $\xi(s;t,x) = (\xi_1(s;t,x),... ,\xi_n(s;t,x))$, where the variables after the semicolon indicate the initial conditions $x= (x_1, ..., x_n)$ at $t$.
We define the smooth characteristic forward flow
$$ 
    \chi: [0,T] \times \mathbb{R}^n \rightarrow \mathbb{R}^n. 
    \quad  (s,x) \mapsto \xi(s;0,x)
$$
Note that $\chi$ satisfies the relation ($d_x$ denoting the Jacobian with respect to the $x$ variables)
$$
  \d_t \chi(t,x) =  d_x \chi(t,x) \cdot a(t,x) \qquad 
    \forall (t,x) \in [0,T] \times \R^n,
$$
which follows upon differentiation of the characteristic differential equations and the initial data with respect to $t$ and $x_k$ ($k=1,\ldots,n$).
Using this relation a straightforward calculation shows that the distributional solution $u  \in  C^\infty([0,T];\mathcal{D}'(\mathbb{R}^{n}))$ to 
$$
  L u  = 0 , \quad u(0) = u_0 \in\mathcal{D}'(\mathbb{R}^n)
$$
is given by
$$
  \dis{u(t)}{\psi} := \dis{u_0}{\psi(\chi(t,.))} \qquad 
    \forall \psi\in \D(\R^n), 0 \leq t \leq T.
$$
If there is a further zero order term $b \cdot u$ in the differential operator $L$, then the above solution formula is modified by an additional factor involving $b$ and $\chi$ accordingly.

In a physical interpretation the characteristic curves correspond to the trajectories of point particles. This provides an idea for introducing a generalized solution concept when the partial differential operator has non-smooth coefficients:
As long as a continuous flow can be defined, the right-hand side in the above definition of $u$ is still meaningful when we 
assume $u_0 \in \mathcal{D}^{\prime 0}(\mathbb{R}^n)$. The distribution $u$ defined in such a way belongs to $AC([0,T]; \mathcal{D}^{'0}(\mathbb{R}^{n}))$ and will be called a \emph{measure solution}.

This approach is not limited to classical solutions of the characteristic system of ordinary differential equations, but can be extended to  more general solution concepts in ODE theory (for example, solutions in the sense of Filippov). Although such a generalized solution will lose the property of solving the partial differential equation in a distributional sense it is a useful generalization with regard to the physical picture.

\subsection{Caratheodory theory}
Let $T>0$ and $\Omega_T = ]0,T[ \times \R^n$. 
Classical Caratheodory theory (cf.\ \cite{Filippov:88}[Chapter 1]) requires the coefficient $a=(a_1, ... , a_n)$ to satisfy
\begin{enumerate}
\item $a(t,x)$ is continuous in $x$ for almost all $t \in [0,T]$,
\item $a(t,x)$ is measurable in $t$ for all fixed $x\in \mathbb{R}^n$ and
\item $\sup_{x\in\mathbb{R}^n} |a(t,x)| \le \beta(t)$ almost everywhere for some positive function $\beta \in L^1([0,T])$.
\end{enumerate}
Then the existence of an absolutely continuous characteristic curve $\xi=(\xi_1, ..., \xi_n)$, which fulfills the ODE almost everywhere, is guaranteed.
Note that the first two Caratheodory conditions ensure Lebesgue measurability of the composition $s\mapsto a(s,f(s))$ for all $f\in {AC}([0,T])^n$, while the
third condition is crucial in the existence proof.\\
A sufficient condition for forward uniqueness of the characteristic system 
is the existence of a positive $\alpha \in L^1([0,T])$,
such that ($\langle . , .\rangle$ denoting the standard inner product on $\mathbb{R}^n$)
\begin{eqnarray*}
  \langle a(t,x) - a(t,y) , x-y \rangle \le \alpha(t) |x-y|^2 
\end{eqnarray*}
for almost all $(t,x), (t,y) \in \ovl{\Omega_T}$ (cf.\ \cite[Theorem 3.2.2]{AgarLak:93}). As well-known from classical ODE theory, forward uniqueness of the characteristic curves yields a continuous forward flow
$$
  \chi: [0,T] \times \mathbb{R}^n \rightarrow \mathbb{R}^n.
  \quad (s,x) \mapsto  \xi(s;0,x)
$$
It is a proper map and for fixed time $\chi(t,.)$ is onto. For the sake of simplicity we assume $a \in C([0,T] \times \mathbb{R}^n)^n$
and $b \in C([0,T] \times \mathbb{R}^n)$.

Let
$$
h_b(t,x) := \exp{ \left( -\int_0^t  b(\tau,\chi(\tau,x)) \,d\tau \right) },
$$
then $u\in \mathcal{D}^{'}(\Omega_T)$ defined by 
\begin{eqnarray} \label{carflowsol}
\langle u, \varphi \rangle_{\mathcal{D}^{'}(\Omega_T)}  := \int_0^T \langle u_0 , \varphi(t,\chi(t,\cdot)) h_b(t,\cdot) \rangle_{\mathcal{D}^{'0}(\mathbb{R}^n)} dt
\end{eqnarray}
(note that $u$ can be regarded as element in $AC([0,T];\mathcal{D}^{'0}(\mathbb{R}^n))$, so the restriction $u(0)$ is well-defined and equal to $u_0\in\mathcal{D}^{'0}(\mathbb{R}^n)$)
solves the initial value problem 
$$
  Lu:= \partial_t u + 
  \sum_{k=1}^n  \partial_{x_k} (a_k \cdot u) +b u = 0, \ \ u(0)=u_0 
$$
on $\Omega_T$, where $a_k \cdot u$ and $b \cdot u$ denotes the distributional product defined by
\begin{eqnarray*}
\cdot: C(\Omega_T)  \times \mathcal{D}^{'0}(\Omega_T) &\rightarrow& \mathcal{D}^{'0}(\Omega_T)\\
  		(f,u) & \mapsto & (\varphi \mapsto \langle u, f \cdot \varphi \rangle_{\mathcal{D}^{'0}(\Omega_T)}).  
\end{eqnarray*}

Applying $L$ on $u$ we obtain
\begin{multline*}
\langle L u, \varphi \rangle_{\mathcal{D}^{'}(\Omega_T)} = \langle u, -\partial_t \varphi - \sum_{k=1}^n a_k \partial_{x_k} \varphi + b \varphi \rangle_{\mathcal{D}^{'}(\Omega_T)} \\
=
 \int_0^T \langle u_0 , (-\partial_t \varphi - \sum_{k=1}^n a_k \partial_{x_k} \varphi + b \varphi)(t,\chi(t,\cdot)) h_b(t,\cdot) \rangle_{\mathcal{D}^{'0}(\mathbb{R}^n)} dt.
\end{multline*}
Set $\phi(t,x):= \varphi(t,\chi(t,x))$ and $\psi(t,x) := \phi(t,x) \cdot h_b(t,x)$, then we have
$$
\partial_t \phi(t,x) = \partial_t \varphi(t,\chi(t,x)) = (\partial_t \varphi + \sum_{k=1}^n a_k(t,x) \partial_{x_k} \varphi)(t,\chi(t,x)),
$$ and
\begin{multline*}
\partial_t \psi(t,x)  = \partial_t \phi(t,x) h_b(t,x) + \phi(t,x) \partial_t h_b(t,x)  \\
= 
(\partial_t \varphi + \sum_{k=1}^n a_k(t,x) \partial_{x_k} \varphi)(t,\chi(t,x)) \cdot  h_b(t,x)- \varphi(t,\chi(t,x)) b(t,\chi(t,x)) h_b(t,x) \\
=(\partial_t \varphi + \sum_{k=1}^n a_k(t,x) \partial_{x_k} \varphi - b \varphi)(t,\chi(t,x)) \cdot h_b(t,x),  
\end{multline*}
thus
$$
\langle L u, \varphi \rangle_{\mathcal{D}^{'}(\Omega_T)} = -\int_0^T \langle u_0 , \partial_t \psi(t,\cdot) \rangle_{\mathcal{D}^{'0}(\mathbb{R}^n)} dt = 
-\int_0^T \partial_t \langle u_0 ,  \psi(t,\cdot) \rangle_{\mathcal{D}^{'0}(\mathbb{R}^n)} dt =0.
$$
for all $\varphi \in \mathcal{D}(\Omega_T)$. The initial condition $u(0)=u_0$ is satisfied, since $\chi(0,x)=x$, thus $h_b(0,x)=1$.

\begin{remark}
In this sense, we can obtain a distributional solution for the Cauchy problem 
\begin{eqnarray*}
P v:= \partial_t v +\sum_{k=1}^n a_k \partial_{x_k} v + c v = 0, \ \ v(0)=v_0,
\end{eqnarray*}
whenever $a \in C([0,T]\times \mathbb{R}^n)^n$ and $c \in \mathcal{D}'([0,T] \times\mathbb{R}^n)$, such that $-{\rm div}(a) + c \in C([0,T] \times \mathbb{R}^n)$ and
$v_0 \in \mathcal{D}^{'0}(\mathbb{R}^n)$. We simply set $b:= -{\rm div}(a) + c$ and construct the solution as above. In other words, such a solution
solves the equation in a generalized sense, relying on the definition of the action of $Q:=\sum_{k=1}^n a_k \partial_k +c$ on a distribution of order $0$ by
\begin{eqnarray*}
\langle Q v, \varphi \rangle_{\mathcal{D}^{'}(\Omega_T)} := -\langle v , \sum_{k=1}^n a_k \partial_{x_k} \varphi \rangle_{\mathcal{D}^{'0}(\Omega_T)} - \langle v, (-{\rm div}(a) +c) \varphi \rangle_{\mathcal{D}^{'0}(\Omega_T)}. 
\end{eqnarray*}
In case where ${\rm div}(a)$ and $c$ are both continuous, we can define the operator $Q$ classically by using the product $\cdot:
\mathcal{D}^{'0}(\Omega_T) \times C(\Omega_T) \rightarrow \mathcal{D}^{'0}(\Omega_T) $ as above.
\end{remark}

\subsection{Filippov generalized characteristic flow}

As we have seen in the previous subsection, forward unique characteristics give rise to a continuous forward flow. But
in order to solve the characteristic differential equation in the sense of Caratheodory, we needed continuity of the coefficient $a$ in the space variables for almost all $t$. 
In case of more general coefficients $a\in L^1_{\rm loc}(\mathbb{R}, L^{\infty}(\mathbb{R}^n))^n$ we can employ
the notion of Filippov characteristics, which replaces the ordinary system of differential equations by a system of differential inclusions (cf. \cite{Filippov:88}). 
The generalized solutions are still absolutely continuous functions.
Again, the forward-uniqueness condition on the coefficient $a$ 
\begin{eqnarray} \label{forward_uniqueness_condition}
\langle a(s,x) - a(s,y) , x-y \rangle \le \alpha(s) |x-y|^2 
\end{eqnarray}
almost everywhere yields unique solutions in the Filippov generalized sense. The generated Filippov flow is again continuous and will
enable us to define measure-valued solutions of the PDE (cf. \cite{PoupaudRascle:97}), as before.

In the Filippov solution concept the coefficient is replaced by a set-valued function $(t,x) \rightarrow A_{t,x} \subseteq \mathbb{R}^n$. It has to have some basic properties which imply the solvability of the resulting system of differential inclusions
\begin{eqnarray*}
\dot{\xi_F}(s) \in A_{s, \xi_F(s)}, \ a.e., \ \ \xi_F(t)=x,
\end{eqnarray*} 
with $\xi_F \in AC([0,\infty[)^n$. These basic conditions are
\begin{enumerate}
\item $A_{t,x}$ is non-empty, closed, bounded and convex for all $x \in \R^n$ and almost all $t\in [0,T]$,
\item $\{ t\in [0,T] \mid \sup_{a\in A_{t,x}}{ \langle a, w \rangle}  < \rho \}$ is Lebesgue measurable for all $x\in \R^n, w\in \R^n, \rho \in \R$, 
\item for almost all $t\in [0,T]$, the set $\bigcup\limits_{x\in K} \{x\} \times A_{t,x}$ is a compact subset of $\R^n \times \R^n$ for $K\csub \R^n$, and
\item there exist a positive function $\beta\in L^1([0,T])$ such that $\sup_{a\in A_{t,x}} |a| \le \beta(t)$ for almost all $t \in[0,T]$ and all $x\in\R^n$.
\end{enumerate}
There are several ways to obtain such a set-valued function $A$ from a coefficient $a\in L^1([0,T]; L^{\infty}(\mathbb{R}^n) )^n$, such that
the classical theory is extended in a compatible way. 
Thus the corresponding set-valued function $A$ should fulfill $A_{t_0,x_0}:=\{a(t_0,x_0)\}$ whenever $a$ is continuous at $(t_0,x_0)\in [0,\infty[ \times \mathbb{R}^n$.

One way obtaining a set-valued function corresponding to a $a \in L^{1}([0,T]; L^{\infty}(\mathbb{R}^n))^{n}$ is
by means of the essential convex hull $\rm{ech}(a)$.
It is defined at $(t,x) \in [0,T]\times\mathbb{R}^{n}$ by
$$
({\rm ech}(a))_{t,x} := \bigcap_{\delta > 0} \bigcap_{
{N\subseteq \mathbb{R}^n \atop \lambda(N)=0}
} {\rm ch}(a(t, B_{\delta}(x) / N)) 
$$
where $\rm ch(M)$ denotes the convex hull of a set $M\subseteq \mathbb{R}^n$ and $\lambda$ is the Lebesgue measure on $\mathbb{R}^n$.

Another way is to use a mollifier $\rho \in \mathcal{S}(\mathbb{R}^n)$ with $\int \rho(x) dx =1$, put $\rho_{\varepsilon}(x) = \varepsilon^{-n} \rho(\varepsilon^{-1} x)$ and
$A_{\varepsilon} := \widetilde{a} \ast \rho_{\varepsilon}\mid_{[0,T] \times \mathbb{R}^n}$ where $\widetilde{a} \in L^{\infty}(\R^{n+1})^n$ is the extension of a function $a\in L^{\infty}([0,T] \times \R^n)^n$ by zero. Then the concept of a generalized graph $C_A$ as defined in \cite{HalSim:07} yields a set-valued function satisfying the above basic properties.

\subsubsection{Measure solutions according to Poupaud-Rascle} \label{section_poupaud_rascle}

Let $\Omega_{\infty} := \,]0,\infty[ \times \mathbb{R}^n$.
We assume $a \in L^1_{\rm loc}(\mathbb{R}_{+}; L^{\infty}(\mathbb{R}^n))^n$ to be a coefficient satisfying the forward uniqueness criterion \eqref{forward_uniqueness_condition}.
Let $L u := \partial_t u + \sum_{i=1}^n \partial_{x_i} (a_i u)$ and 
$\xi_F$ be the unique solution to
\begin{eqnarray}
\dot{\xi_F}(s) & \in & {\rm ech}(a)_{s, \xi_F(s)}, \ \ \xi_F(t) = x.
\end{eqnarray}
 The map
$$
    \chi_F: \mathbb{R}_{+} \times \mathbb{R}^n \rightarrow 
    \mathbb{R}^n, \quad (t,x) \mapsto \xi_F(t;0,x)
$$
is the continuous Filippov (forward) flow.

\begin{definition}[Solution concept according to Poupaud-Rascle]
Let $u_0\in \mathcal{M}_b(\mathbb{R})^n$ be a bounded Borel measure, then the image measure at $t\in[0,\infty[$ is 
\begin{eqnarray} \label{filflowsol}
u(t)(B) :=  \int_{\mathbb{R}^n} 1_{B}(\chi_F(t,x)) du_0(x), 
\end{eqnarray}
where $B\subseteq \mathbb{R}^n$ is some Borel set. The map $u \col [0,\infty[ \to \mathcal{M}_b(\mathbb{R}^n))$ belongs to $C([0,\infty[; \mathcal{M}_b(\mathbb{R}^n))$ and is called
a measure solution in the sense of Poupaud-Rascle of the initial value problem
\begin{eqnarray*} \label{pdo2}
L u&:=& \partial_t u + \sum_{k=1}^n \partial_{x_k} (a_k \cdot u) = 0, \ \ u(0) = u_0.
\end{eqnarray*}
Note that $u$ defines a distribution of order $0$ in $\mathcal{D}'(\Omega_{\infty})$ by
\begin{eqnarray*}
\langle u, \varphi \rangle_{\mathcal{D}'(\Omega_{\infty})} := \int_0^{\infty} \langle u_0,  \varphi(t,\chi_F(t,x)) \rangle_{\mathcal{D}^{'0}(\mathbb{R}^n)}\, dt, \quad \quad \forall \varphi\in \mathcal{D}'(\Omega_{\infty}). 
\end{eqnarray*}
\end{definition}

The solution concept of Poupaud-Rascle does not directly solve the partial differential equation in a distributional sense,
but it still reflects the physical picture of a "transport process" as imposed by
the properties of the Filippov characteristics. Nevertheless, in the cited paper of Poupaud-Rascle(\cite{PoupaudRascle:97}) the authors present an a posteriori definition of the particular product $a\cdot u$, which restore the validity of the PDE in a somewhat artificial way. We investigate this in the sequel in some detail.

\begin{definition}[A posteriori definition of a distributional product in the sense of Poupaud-Rascle]
Let $u\in \mathcal{D}^{\prime}(\mathbb{R}^n)$ be a distribution of order $0$ and $a\in {L^1_{\rm loc}([0,\infty[, L^{\infty}(\mathbb{R}^n))}^n$,
satisfying the forward uniqueness condition (\ref{forward_uniqueness_condition}), such that there exists a continuous Filippov flow $\chi_F$. Furthermore we assume that $u$ is a generalized solution of the initial value problem
as defined in (\ref{filflowsol}). Then we define the product $a \bullet u =(a_k\cdot u)_k$ in $\mathcal{D}'(]0,\infty[\times \mathbb{R}^n)^n$ by
\begin{eqnarray*}
\langle a \bullet u, \varphi \rangle_{\mathcal{D}'(\Omega_{\infty})} &:=& \langle u_0, \int_0^{\infty} \partial_t \chi_F(t,x) \varphi(t, \chi(t,x) ) dt \rangle_{\mathcal{D}^{'0}(\mathbb{R}^n)}, \ \ \varphi \in \mathcal{D}(\Omega_{\infty}).
\end{eqnarray*}
\end{definition}

\begin{remark}
Note that the product $a\cdot u$ is defined only for distributions $u$ that are a generalized solutions (according to Poupaud-Rascle) of the initial value problem 
(\ref{pdo2}) with the coefficient $a$. The domain of the product map $(a,u) \mapsto a\bullet u$, as subspace of $\mathcal{D}^{'0}(\R^n) \times \mathcal{D}^{'0}(\R^n)$  has a complicated structure: Just note that the property to generate a continuous characteristic Filippov flow $\chi_F$
is not conserved when the sign of the coefficient $a$ changes, as we have seen for the coefficient $a(x)={\rm sign}(x)$.
\end{remark}

\begin{example} \label{example_poupaud_rascle_1}
Consider problem (\ref{pdo2}) with the coefficient $a(x):= -{\rm sign}(x)$ subject to the initial condition $u_0= 1$. Then the continuous Filippov flow
 is given by 
$$
 \chi_F(t,x) = -(t+x)_{-} H(-x) + (x-t)_{+} H(x).
$$
We have $\chi_F(t,0)= t_{+}- (-t)_{-} = 0$
and
\begin{eqnarray*}
\partial_t \chi_F(t,x) = - H(-t-x) H(-x) - H(x-t) H(x) \ \ \text{for almost all } \ t\in[0, \infty[. 
\end{eqnarray*}
The generalized solution $u$ is defined by $\langle u, \varphi \rangle: = \int_0^{\infty} \langle u_0, \phi(t, x) \rangle \,dt$,
where $\phi(t,x) := \varphi(t, \chi(t,x))$. We have that
\begin{eqnarray*}
\phi(t,x) := \left\{  \begin{array}{ll}
\varphi(t,x+t) & x\le 0, 0\le t \le -x \\
\varphi(t,0) & t \ge |x| \\
\varphi(t,x-t) & x\ge 0, 0\le t \le x ,\\
\end{array}
\right. 
\end{eqnarray*}
thus 
\begin{multline*}
\langle u, \varphi \rangle_{\mathcal{D}^{'0}(\Omega_\infty)} :=\int_0^{\infty} \langle u_0, \phi(t, x) \rangle \,dt
= \int_0^{\infty} \int_{-\infty}^{\infty}  \phi(t, x) \,dx\ \,dt \\
= 
2 \int_0^{\infty} \varphi(t, 0) t \,dt  +  \int_0^{\infty} \int_{-\infty}^{-t} \varphi(t,x+t) \,dx\ \,dt +
\int_0^{\infty} \int_{t}^{\infty} \varphi(t,x-t)  \,dx\ \,dt \\
= 2 \int_0^{\infty} \varphi(t, 0) t \,dt  +  \int_0^{\infty} \int_{-\infty}^{0} \varphi(t,z) \,dz\ \,dt +
\int_0^{\infty} \int_{0}^{\infty} \varphi(t,z)  \,dz\ \,dt 
= \langle 1+ 2t \delta, \varphi(t,\cdot) \rangle_{\mathcal{D}^{'0}(\mathbb{R}^n)} 
\end{multline*}
This generalized solution gives rise to the following product
\begin{eqnarray*}
\langle (-{\rm sign}(x)) \bullet (1+2t \delta(x)), \varphi \rangle_{\mathcal{D}'(\Omega_{\infty})} &:=&
 \langle 1 , \int_0^\infty  \partial_t \xi_F(t,x) \varphi(t, \xi_F(t,x)) dt \rangle_{\mathcal{D}^{'0}(\mathbb{R}^n)}
\end{eqnarray*}
in $\mathcal{D}'(\Omega_{\infty})$. Evaluating the right-hand side we obtain 
\begin{multline*}
\langle 1 , \int_0^\infty  \partial_t \chi_F(t,x) \varphi(t, \chi_F(t,x)) dt \rangle 
=
\int_{-\infty}^{\infty}  \left( -\int_{0}^{\infty} H(-x)H(-x-t) \varphi(t, x+t) \,dt \right. \\
\left.
-\int_{0}^{\infty} H(x)H(x-t) \varphi(t, x-t) \,dt \right) \,dx.
\end{multline*}
Since $H(-x)H(-x-t) =H(-x-t)$ and $H(x) H(x-t) = H(x-t)$ for $t\ge 0$ the latter gives upon substitution
\begin{eqnarray*}
 \langle 1 , \int_0^\infty  \partial_t \chi_F(t,x) \varphi(t, \chi_F(t,x)) dt \rangle 
&=& -\int_{-\infty}^{\infty} {\rm sign}(z)  \int_{0}^{\infty} \varphi(t, z) \,dt  \,dz,
\end{eqnarray*}
hence
\begin{eqnarray*}
(-{\rm sign}(x)) \bullet (1+2t \delta(x)) = -{\rm sign}(x).
\end{eqnarray*}
However, we cannot define the product if  $-{\rm sign}(x)$ is replaced by $+{\rm sign}(x)$, since
the Filippov characteristics $\xi_F(t;0,x)$ are no longer forward unique and thus do not generate a continuous  Filippov flow $\chi_F$.
\end{example}

\begin{example}
We consider the same coefficient $a(x):=- {\rm sign}(x)$ as before, but now we set $u_0:=\delta$. We obtain the generalized solution
\begin{eqnarray*}
\langle u, \varphi \rangle_{\mathcal{D}^{'}(\Omega_{\infty})} &:=&  {\langle 1 \otimes \delta, \varphi(t, \chi_F(t,x)) \rangle}_{\mathcal{D}^{'}(\Omega_{\infty})} = \int_0^\infty  \varphi(t, \chi_F(t,0))  \,dt 
\end{eqnarray*}
This enables us to calculate the product
\begin{eqnarray*}
&&\langle (- {\rm sign}(x)) \bullet \delta(x),\varphi \rangle = -\langle \delta, \int_0^{\infty} \partial_t \chi_F(t,x) \varphi(t, \xi_F(t,x)) dt \rangle.
\end{eqnarray*}
Putting $\psi(x)=\int_0^{\infty} \partial_t \xi_F(t,x) \varphi(t, \xi_F(t,x)) dt$ and observe that
\begin{eqnarray*}
\psi(x) &:=&\int_0^{\infty} \partial_t \chi_F(t,x)  \varphi(t, \chi_F(t,x)) dt = 
\int_0^{-x} \varphi(t, x+t ) dt , \ \text{if}\ x<0
\end{eqnarray*}
and
\begin{eqnarray*}
\psi(x) &=&\int_0^{\infty} \partial_t \chi_F(t,x)  \varphi(t, \chi_F(t,x)) dt = 
-\int_0^{x}  \varphi(t,x-t) dt,\ \text{if}\ x>0.
\end{eqnarray*}
At $x=0$ we obtain $\psi(0)=\lim_{x\rightarrow 0_{-}} \psi(x)=\lim_{x\rightarrow 0_{+}} \psi(x)=0$, so it follows that
$(- {\rm sign}) \bullet \delta = 0$.
\end{example}

\begin{example}\label{test1}
Let $a(t,x):= 2 H(-x)$, so that the Filippov flow is given by
\begin{eqnarray*}
\chi_F(t,x) = -(x+2 t)_{-} H(-x) + x H(x). 
\end{eqnarray*}
We have $\chi_F(t,0)= -2t_{-} =0$ and 
\begin{eqnarray*}
\partial_t \chi_F(t,x):=2 H(-x-2t) H(-x).
\end{eqnarray*}
Hence $\partial_t \chi_F(t,0) = 0$ for almost all $t \in [0, \infty[$.
If $u_0=1$ the generalized solution is
\begin{eqnarray*}
\langle u, \varphi \rangle_{\mathcal{D}^{'0}(\Omega_{\infty})} := \int_0^{\infty} \int_{-\infty}^{\infty} \phi(t, x) \,dx \,dt, 
\end{eqnarray*}
where $\phi(t,x) = \varphi(t, \chi_F(t,x))$. Since
\begin{eqnarray*}
\phi(t,x)\mid_{\{x< -2t\}} &=& \varphi(t, x+2t) \\
\phi(t,x)\mid_{\{-2t \le x \le 0 \}} &=& \varphi(t, 0) \\
\phi(t,x)\mid_{\{0 < x   \}} &=& \varphi(t, x) ,
\end{eqnarray*}
we obtain
\begin{multline*}
\langle u, \varphi \rangle_{\mathcal{D}^{'0}(\Omega_{\infty})} 
= \int_0^{\infty}
 \left( \int_{-\infty}^{-2t} \varphi(t, x+2t) \,dx  + 2t \varphi(t, 0)  + \int^{\infty}_{0} \varphi(t, x) \,dx \right) \,dt 
= \langle 1+ t \delta, \varphi(t,\cdot) \rangle_{\mathcal{D}^{'0}(\mathbb{R}^n)},
\end{multline*}
hence $u= 1+2 t \delta(x)$.
Again we determine the product $(2H(-x)) \bullet  (1+2 t \delta(x))$ by
\begin{multline*}
\langle 2 H(-x) \bullet  (1+2t \delta(x)), \varphi \rangle_{\mathcal{D}'(\Omega_{\infty})} = 
2 \int_{-\infty}^{\infty} \int_0^{\infty} H(-x) H(-x-2t) \varphi(t, x+2t) \,dt \,dx \\
=
2 \int_{-\infty}^{\infty} \int_0^{\infty}  H(-x-2t) \varphi(t, x+2t) \,dt \,dx 
= 2 \int_0^{\infty} \int_{-\infty}^{\infty} H(-z) \varphi(t,z)  \,dz \,dt = \langle 1\otimes 2 H(-\cdot), \varphi \rangle_{\mathcal{D}'(\Omega_{\infty})}.
\end{multline*}
We obtain $(2H(-x)) \bullet  (1+2 t \delta(x)) = 2 H(-x)$. 
Observe that together with the result in Example (\ref{example_poupaud_rascle_1}) $(-{\rm sign}(x)) \bullet (1+2 t \delta(x))=(2 H(-x) - 1 ) \bullet (1+2 t \delta(x))$ we 
can conclude that either $(-1) \bullet  (1+2 t \delta(x)) $ is not defined or the product $\bullet$ is not distributive.
In fact , it is not difficult to see that $(-1) \bullet  (1+2 t \delta(x)) $ cannot be defined in this way, neither can $1 \bullet  (1+2 t \delta(x))$.
\end{example}

\begin{example}[generalization of Example \ref{example_poupaud_rascle_1}] \label{poupaud_rascle_jump} 
Let $c_1 \ge c_2$ be two constants, and $\alpha \in [c_1, c_2]$.
Consider the $a(t,x):=c_1 H(\alpha t- x) + c_2 H(x- \alpha t)$. We set $t_1(x):= \frac{-x}{c_1-\alpha}$ if $x<0$ and $t_2(x):= \frac{x}{\alpha-c_2}$  for $x>0$. 
The unique Filippov flow is given by
\begin{eqnarray*}
\chi_F(t,x) =\left\{ 
\begin{array}{ll}
c_1 t + x  &  x < 0, t < t_1(x)\\
\alpha t & x<0, t\ge t_1(x) \\
\alpha t & x=0, \\
c_2 t +x & x>0, t\le t_2(x) \\
 \alpha t & x>0 , t\ge t_2(x)
\end{array}\right.
\end{eqnarray*}
The generalized solution of the initial value problem $Lu:=\partial_t u + \partial_x (a \cdot u)=0, \ \ u(0)=u_0\in L^1_{\rm loc}(\mathbb{R})$,
according to Poupaud-Rascle is given by 
\begin{multline*}
\langle u, \varphi\rangle_{\mathcal{D}^{'0}(\Omega_T)}=
\int_0^T \langle u_0, \varphi(t, \chi_F(t,\cdot)) \rangle_{\mathcal{D}^{'0}(\mathbb{R})} \,dt = \int_{-\infty}^0 \int_{0}^{t_1(x)} u_0(x)\varphi(t,c_1 t+x ) \,dt\ \,dx \\
+
\int_{-\infty}^0 \int_{t_1(x)}^{T} u_0(x)\varphi(t,\alpha t ) \,dt\ \,dx + 
\int_{0}^{\infty} \int_{0}^{t_2(x)} u_0(x)\varphi(t,c_2 t+x ) \,dt\ \,dx 
+
\int_{0}^{\infty} \int_{t_2(x)}^{T} u_0(x) \varphi(t,\alpha t ) \,dt\ \,dx  \\
=
  \int_{0}^{T}\int_{-\infty}^{-t(c_1-\alpha)} u_0(x)\varphi(t,c_1 t+x )  \,dx\ \,dt 
+
\int_{0}^{T}  \int_{-t(c_2 -\alpha)}^{\infty} u_0(x)\varphi(t,c_2 t+x )  \,dx\ \,dt \\ +  \int_{0}^T \left(\int_{-t(c_1-\alpha)}^{t(\alpha-c_2)} u_0(x)  \,dx \right) \varphi(t,\alpha t ) \,dt,
\end{multline*}
hence 
\begin{eqnarray*}
u:= u_0(x-c_1 t) H(\alpha t-x) + u_0(x-c_2 t) H(x-\alpha t) + \left(\int_{-t(c_1-\alpha)}^{t(\alpha-c_2)} u_0(x) \,d x\right) \delta(x-\alpha t).
\end{eqnarray*}
\end{example}

\subsection{Colombeau generalized flow}

In this subsection we consider the solvability of the ordinary differential equations
for the characteristics in the setting of Colombeau generalized functions.
Our main focus will be on distributional shadows of such generalized solutions. It will appear that
under certain assumptions on the right-hand side, the distributional shadow exists and is absolutely continuous. We will also 
show a uniqueness result for distributional shadows.

\begin{theorem}[Existence]\label{colombeau_existence}
Assume $A  \in \G(\overline{\Omega_T})^n$ with a representative $(A_{\varepsilon})_{\varepsilon}$, such that
\begin{eqnarray} \label{colombeau_existence_bound}
\sup_{x\in\mathbb{R}^n}|A_{\varepsilon}(t,x)| \le  \beta(t) ,\ \ \varepsilon \in\, ]0,1], \ \text{almost everywhere in}\ t\in[0,T]
\end{eqnarray}
holds, where $\beta$ is some positive function in $L^1([0,T])$. Let $(\widetilde{t},\widetilde{x}) \in \widetilde{\ovl{\Omega_T}}$
be the initial condition.
Then there exists a c-bounded solution $\xi \in \G([0,T])^n$ to the initial value problem
\begin{eqnarray*}
\dot{\xi}(s) &=& A(s,\xi(s)), \ \xi(\widetilde{t})=\widetilde{x}.
\end{eqnarray*}
Furthermore, there exists some $(t,x) \in \ovl{\Omega_T}$, $\xi_C \in AC([0,T])$ such that for any representantive $(\xi_{\varepsilon})_{\varepsilon}$ of $\xi$, $(t_\varepsilon, x_{\varepsilon})_{\varepsilon}$ of $(\widetilde{t},\widetilde{x})$ there exists subsequences $(t_{\varepsilon_j}$, $x_{\varepsilon_j})_j$,$(\xi_{\varepsilon_j})_j$ with 
$\lim_{j \rightarrow \infty} (t_{\varepsilon_j}, x_{\varepsilon_j}) =(t,x)$ and $\xi_{\varepsilon_j} \stackrel{j\rightarrow \infty}{\rightarrow} \xi_C$ uniformly on $[0,T]$ and $\xi_C(t)=x$. 
\end{theorem}
\begin{proof}
By classical existence and uniqueness we obtain $\xi_{\varepsilon}$ for each $\varepsilon \in ]0,1]$ such that
\begin{eqnarray*}
\xi_{\varepsilon}(s) &=& x_{\varepsilon} + \int_{t_{\varepsilon}}^s A_{\varepsilon}(\tau,\xi_{\varepsilon}(\tau))\,d\tau 
\end{eqnarray*}
holds. Condition (\ref{colombeau_existence_bound}) yields $|\xi_{\varepsilon}(s)| \le |x_{\varepsilon}| + |\int_t^s \beta(\tau) \,d\tau|$ for all
$s,t\in[0,T]$, hence c-boundedness of $(\xi_{\varepsilon})_{\varepsilon}$ on $[0,T]$ and furthermore moderateness of $\dot{\xi_{\varepsilon}}$ (by \cite{GKOS:01}[Proposition 1.2.8]). In fact this existence result is quite similar to the one given in \cite{GKOS:01}[Proposition 1.5.7].

To prove the existence of a convergent subsequence of  $(\xi_{\varepsilon})_{\varepsilon}$, we may assume 
without loss of generality that $\lim_{\varepsilon \rightarrow 0}(t_{\varepsilon}, x_{\varepsilon})=(t,x) \in \overline{\Omega_T}$. Note that the family $(\xi_{\varepsilon})_{\varepsilon}$
is uniformly bounded and equicontinuous, since
\begin{eqnarray*}
|\xi_{\varepsilon}(s) - \xi_{\varepsilon}(s')| \le  |\int_s^{s'} \beta(\tau) \,d\tau| \quad s,s' \in [0,T], \varepsilon \in ]0,1].
\end{eqnarray*}
The Theorem of Arzela-Ascoli yields a subsequence $(\xi_{\varepsilon_j})_{j}$ converging uniformly to some $\xi_C \in C([0,T])$. Clearly,
$\xi_C(t)=\lim_{j\rightarrow \infty} \xi_{\varepsilon_j}(t_{\varepsilon_j})=\lim_{j\rightarrow \infty} x_{\varepsilon_j}=x$. 
such that
\begin{eqnarray*}
\lim_{j\rightarrow \infty} \sup_{s \in [0,T]} |\xi_{\varepsilon_j}(s) -\xi_C(s)| = 0,
\end{eqnarray*}
holds. We have for all $s,s'\in[0,T]$,
\begin{multline*}
|\xi_C(s) - \xi_C(s')| \le |\xi_C(s) - \xi_{\varepsilon_{j}}(s)|+ |\xi_{\varepsilon_{j}}(s)-\xi_{\varepsilon_{j}}(s')| +|\xi_{\varepsilon_{j}}(s') - \xi_C(s')|  \\
\le |\xi_C(s) - \xi_{\varepsilon_{j}}(s)| +\int_{s}^{s'} \beta(\tau) \,d\tau +|\xi_{\varepsilon_{j}}(s') - \xi_C(s')| \stackrel{j \rightarrow \infty}{\rightarrow}\int_{s}^{s'} \beta(\tau) \,d\tau , 
\end{multline*}
hence $\xi_C$ is absolutely continuous on $[0,T]$.
\end{proof}

\subsection{Semigroups defined by characteristic flows}

Let $X$ be a Banach space and $(\Sigma_t)_{t\in [0,\infty[}$ be a family of bounded operators $\Sigma_t$ on $X$.
Consider the following conditions: 
\begin{enumerate}
\item $\Sigma_0 = {\rm id}$ 
\item $\Sigma_s \circ \Sigma_t = \Sigma_{s+t}$ for all $s,t \in [0,\infty[$ and 
\item the orbit maps 
\begin{eqnarray*}
\sigma_{u_0}: [0,\infty[ &\rightarrow& X \\ 
t &\mapsto& \Sigma_t(u_0)
\end{eqnarray*}
are continuous for every $u_0 \in X$. 
\end{enumerate}
If (i) and (ii) are satisfied, then we call $(\Sigma_t)$ a semi-group acting on $X$.
If  in addition property (iii) holds, we say  $(\Sigma_t)_{t\in [0,\infty[}$ is a semi-group of type $C_0$.

We briefly investigate how the solution concepts discussed in subsections 1.1 and 1.2 fit into the picture
of semi-group theory when the coefficient $a$ is time-independent. First we return to the classical Caratheodory case: 
Let $a\in C(\mathbb{R}^n)^n$ and assume that $a$ suffices the forward uniqueness condition (\ref{forward_uniqueness_condition}).
This implies that the characteristic flow $\chi : \overline{\Omega_{\infty}} \rightarrow \mathbb{R}^n$ is continuous and $\chi(t,\cdot)$ is onto $\mathbb{R}^n$ for fixed $t\in[0,\infty[$. Furthermore we have $\chi(s,\chi(r,x)) = \chi(s+r,x)$ for all $ x \in \mathbb{R}^n$ and $r,s \in [0,T] $ with $s+r \in [0,T]$, since $a$ is time independent.

Consider the inital value problem $P u = \partial_t + \sum_{k=1}^n a_k \partial_{x_k} u =0$ with initial condition
$u(0) = u_0\in C_{0}(\mathbb{R}^n)$ (i.e. vanishes at infinity). It is easy to verify that
\begin{eqnarray*}
\Sigma_t : C_0(\mathbb{R}^n) &\rightarrow& C_0(\mathbb{R}^n)\\
           u_0 &\mapsto& \chi^{\ast} u_0 
\end{eqnarray*}
defines $C_0$ semigroup on the Banach space $C_0(\overline{\Omega_{\infty}})$:
Note that $\Sigma_t$ is a bounded operator on $C_0(\mathbb{R}^n)$ for each $t\in [0,\infty[$, as $\chi(t,\mathbb{R}^n) = \mathbb{R}^n$, so
\begin{eqnarray*}
\| \Sigma_t (u_0) \|_{\infty} &=& \sup_{x\in \mathbb{R}^n} \|u_0(\chi(t,x)) \| = \sup_{x\in \mathbb{R}^n} \|u_0(\chi(t,x)) \| = \sup_{x\in \mathbb{R}^n} \|u_0(x)\| =\|u_0\|_{\infty}.
\end{eqnarray*}
We have that $\| \Sigma_t  \| = 1$ for all $t\in[0,\infty[$.
Condition (i) and (ii) follow directly from the flow properties of $\chi$. The continuity condition (iii), which is equivalent to
\begin{eqnarray*}
\lim_{t  \rightarrow 0_{+}}\| \Sigma_t(u_0) - u_0\|_{\infty} = \lim_{t  \rightarrow 0_{+}} \sup_{x\in \mathbb{R}^n} \|u_0(\chi(t,x)) - u_0(x)\| = 0,
\end{eqnarray*}
holds, since $(\chi(t,x))_{x\in\mathbb{R}^n}$ is an equicontinuous family and $u_0$ vanishes at infinity.

\begin{remark}
For a coefficient $a$ in $L^{\infty}(\mathbb{R}^n)^n$ we can also define a semi-group on $C_0(\mathbb{R}^n)$  
by $\Sigma_t(u_0):= u_0(\chi_F(t,x))$, where $\chi_F$ is the generalized Filippov flow as introduced earlier. This is due to the fact,
that the Filippov flow has almost the same properties as the Caratheodory flow. 
\end{remark}

It seems natural to understand the  solution concepts as defined by (\ref{carflowsol}) and (\ref{filflowsol})
as action of the dual semigroup $(\Sigma_t^{\ast})$ on 
\begin{eqnarray*}
 \langle u(t), \varphi \rangle_{\mathcal{D}^{'0}(\mathbb{R}^n)} = \langle \Sigma_t^{\ast} u_0, \varphi \rangle_{\mathcal{D}^{'0}(\mathbb{R}^n)} 
=\langle u_0, \Sigma_t (\varphi) \rangle_{\mathcal{D}^{'0}(\mathbb{R}^n)}. 
\end{eqnarray*}
 the Banach space of finite complex Radon measures, the dual space of $C_0(\mathbb{R}^n)$ 
(cf. \cite[Chapter 4]{Davies:80}, \cite[Chapter 1.10]{Pazy:83} or \cite[Chapter IX.13]{Yosida:80} for the general setting).
However, in general the dual semi-group is not of class $C_0$ (cf. \cite[Example 1.31]{Davies:80}). This is only guaranteed if we start from a $C_0$ semi-group  defined on a reflexive Banach space. 

Nevertheless the solution concepts in (\ref{carflowsol}) and (\ref{filflowsol}) still yield the semi group properties (i) and (ii) with 
weak-$\ast$ continuity replacing the strong continuity property (iii).

The situation is much easier with Hilbert spaces, of course. We conclude with an example involving a discontinuous coefficient.


\begin{example} \label{semigroup_L2_example}
Let $a\in L^{\infty}(\mathbb{R})$ such that there exist $c_0,c_1>0$ such that
$c_1 < a(x) < c_2$ almost everywhere. We want to solve the initial value problem
\begin{eqnarray*}
Pu=\partial_t u + a(x) \partial_x u =0, \ \ u(0)=u_0 \in L^2(\mathbb{R})
\end{eqnarray*}
for $u \in AC([0,T]; L^2(\mathbb{R})) \cap L^1([0,T]; H^1(\mathbb{R}))$.

Let $A(x) = \int_{0}^{x} a(y)^{-1} \,dy$, which is Lipschitz continuous and strictly increasing (thus globally invertible) 
and observe that $\chi(t,x) = A^{-1}(t+A(x))$ defines the (forward)characteristic flow that solves
$$
\chi(t,x) = x + \int_0^{t} a(\chi(\tau,x)) \,d\tau.
$$
Let $Q:= - a(x) \partial_x$ with domain $D(Q):=  H^1(\mathbb{R})$.
The resolvent of $Q$ for $\Re(\mu) > 0$ is obtained from the equation
\begin{eqnarray*}
(-Q + \mu ) v = f, \ \ f\in L^2(\mathbb{R}).
\end{eqnarray*}
Upon division by  $a$ we deduce
\begin{eqnarray} \label{characteristic_equation}
\partial_x v +\frac{\mu}{a}{v} = \frac{f}{a}.
\end{eqnarray}

Let us first consider uniqueness: Let  $w \in H^1(\R)$ satisfy
\begin{equation}  \label{semigroup_L2_homogenous}
\partial_x w + \frac{\mu}{a} w = 0.
\end{equation}
Since $w$ is absolutely continuous we have
$$
w(x) = C \exp{\left(- 2 {\rm Re}(\mu) \int_{-\infty}^{x} \frac{1}{a(z) }  \,dz \right)}
$$
for some constant $C$. But $w \in L^2(\mathbb{R})$ if and only if $C=0$, thus $w = 0$.
%
%
%

Existence:
One easily verifies that
\begin{multline*}
v(x)=(R(\mu) f)(x) := \int_{-\infty}^{x} \exp{\left(- \mu \int_{y}^x a(z)^{-1} dz \right)} \frac{f(y)}{a(y)} \, dy  
= \int_{-\infty}^{x} \exp{(- \mu (A(x)-A(y)) )} \frac{f(y)}{a(y)} \, dy    
\end{multline*} 
is a solution of (\ref{characteristic_equation}) in $AC(\mathbb{R})$. 
Upon substitution $y \mapsto z = A(x) - A(y)$ in the right-most integral we obtain
\begin{eqnarray} \label{laplace_transform}
(R(\mu) f)(x) = \int_{0}^{\infty} \exp{(-\mu z)} f(\chi(-z,x))  \,dz,
\end{eqnarray}
which is the Laplace transform of $f(\chi(-.,x))$.

We denote the kernel of the integral operator $R(\mu)$ by
\begin{eqnarray*}
M(x,y) &:=& H(x-y) \exp{(- \mu (A(x)-A(y))} a(y)^{-1}.
\end{eqnarray*} 
We briefly sketch the derivation of $L^2$ estimates for the operator powers $R(\mu)^k$ for ${\rm Re}(\mu)>0$:
Note that $R(\mu)^k$ is an iterated integral operator of the form
\begin{multline*}
 R(\mu)^k f (x) :=R(\mu)^{k-1} \left(\int_{\mathbb{R}} M(\cdot,z_1) f(z_1) \,dz_1 \right) (x)\\= 
 \int_{\mathbb{R}} \cdots  \int_{\mathbb{R}} M(x, z_{k}) M(z_k, z_{k-1}) \ldots M(z_2,z_1) f(z_1) \,dz_{k-1} \ldots dz_2 \, dz_1. 
\end{multline*}
To simplify notation let $z=(z_1, \dots, z_k)$, $d^kz = dz_1 \dots dz_k$, $h(z):=\exp{ (-\mu \sum_{l=1}^k z_l)}$, and 
$g(x,z): = f(\chi(-\sum_{l=1}^k z_l, x))$. Using the flow property of $\chi$ we obtain that
\begin{eqnarray*}
R(\mu)^k f (x) = \int_{[0, \infty[^k} h(z) g(z,x) \, d^k z 
\end{eqnarray*}
holds, hence by by the integral Minkowski inequality
\begin{multline} \label{semigroup_L2_existence}
  \|R(\mu)^k f \|_{L^2} \leq \left( 
  \int_{\mathbb{R}}\left(\int_{[0, \infty[^k} |h(z)| |g(z,x)| \,d^k z \right)^2  \,dx \right)^{1/2}\\
  \leq 
 \int_{[0, \infty[^k}  |h(z)|     
\left(\int_{\mathbb{R}}  |g(z,x)|^2  \,dx \right)^{1/2} \,d^k z.
\end{multline}
Since
$$
  \int_{\mathbb{R}}  |g(z,x)|^2 \,dx  
  = \int_{\mathbb{R}} |f(\chi (-\sum_{l=1}^k z_l, x) )|^2  \,dx
    = \int_{\mathbb{R}} |f(y)|^2  
    \left|\frac{a( y) }{a(\chi (\sum_{l=1}^k z_l, y ))}\right| \,dy 
\le 
\frac{c_1}{c_0} \, \| f\|_{L^2}^2,
$$
we conclude
$$
  \| R(\mu)^k f \|_{L^2} \le  
  \sqrt{\frac{c_1}{c_0}}\cdot \| f\|_{L^2} \cdot
  \int_{]-\infty,0]^k}  |h(z)|    
 \,d^kz = 
   \frac{\sqrt{\frac{c_1}{c_0}}}{{\rm Re(\mu)}^{k}}\cdot \| f\|_{L^2} . 
$$
The  Hille-Yosida theorem 
(\cite[Theorem 5.2]{Pazy:83})
yields that $Q$ generates the $C_0$ semigroup 
\begin{eqnarray*}
 \Sigma_t : L^2(\mathbb{R}) &\rightarrow& L^2(\mathbb{R})\\
     	u_0 &\mapsto& u_0(\chi(-t,x)).
\end{eqnarray*}
The resolvent operator $\mu \mapsto R(\mu)$ (defined for ${\rm Re}(\mu)> 0$) is the Laplace transform of the semigroup $t \rightarrow  \Sigma_t$
as indicated in (\ref{laplace_transform}). 

\end{example}

\begin{remark}
Since $L^2(\mathbb{R}^n)$ is reflexive the dual semigroup is $C_0$ as well and has as its generator the adjoint operator $Q^*$.
\end{remark}

\begin{remark}
If we assume additional regularity on the coefficient, e.g.\ $a\in C^{\sigma}_{\ast}(\mathbb{R})$ with $\sigma>0$,  in Example \ref{semigroup_L2_example}, then we obtain a $C_0$ semigroup $(\Sigma_t)$ acting on the
Hilbert space $H^s(\mathbb{R})$ with $0 \leq s < \sigma$. We may then use the fact that the (square of the) Sobolev norm  $\norm{v}{s}^2$  is equivalent to the following expression  (cf.\ \cite[Equation (7.9.4)]{Hoermander:V1})
$$
   \int |v(x)|^2\,dx + C_s \int \int \frac{|v(x) - v(y)|^2}{|x-y|^{n + 2 s}} \, dx \, dy,
$$
where the constant $C_s$ depends only on the dimension $n$ and $s$. 
From this it can be shown that we may have $D(Q)=H^{s+1}(\mathbb{R})$ as domain of $Q$ (this also corresponds to the special case of the mapping properties stated in \cite[Chapter 2]{Taylor:91}).
Clearly,  uniqueness in the characteristic equation (\ref{semigroup_L2_homogenous}) is still valid. A corresponding variant of the estimate (\ref{semigroup_L2_existence}) for the powers of the resolvent operator $R(\mu)^k$ on $H^s(\mathbb{R})$ is obtained by  the following calculation (with the notation $h$ and $g$ as in Example \ref{semigroup_L2_example}):
\begin{multline*}
 \|R(\mu)^k f \|_{s}^2 = \|R(\mu)^k f \|_{0}^2 +   
  \int_{\mathbb{R}} \int_\mathbb{R} \frac{|(R(\mu)^k f)(x)-(R(\mu)^k f)(y)|^2}{|x-y|^{1+2s}} \,dx \,dy \\
  \le   
 \frac{\frac{c_1}{c_0}}{{\rm Re(\mu)}^{2k}}\cdot \| f\|_{0}^2  +
\int_{\mathbb{R}} \int_\mathbb{R} \left(\, \int_{[0,\infty[^k} |h(z)| \frac{|g(x,z) - g(y,z)|}{|x-y|^{\frac{1+2s}{2}}} \,dz \right)^2  \,dx \,dy 
\\
  \le \frac{\frac{c_1}{c_0}}{{\rm Re(\mu)}^{2k}}\cdot \| f\|_{0}^2
    + \left( \, 
\int_{[0,\infty[^k}  |h(z)| \left( \int_{\mathbb{R}} \int_\mathbb{R}  \frac{|g(x,z) - g(y,z)|^2}{|x-y|^{1+2s}} \,dx \,dy   \right)^{1/2}   \,d^k z \right)^2.
\end{multline*}
To carry out the $x$ and $y$ integrations we use the substitutions $x' = \chi(-\sum_{l=1}^k z_l, x)$, $y' = \chi(-\sum_{l=1}^k z_l, y)$ to obtain
\begin{multline*}
\int_{\mathbb{R}} \int_\mathbb{R}  \frac{|g(x,z) - g(y,z)|^2}{|x-y|^{1+2s}} \,dx \,dy
 =
\int_{\mathbb{R}} \int_\mathbb{R}  \frac{| f( \chi(-\sum_{l=1}^k z_l, x) )-   f(\chi(-\sum_{l=1}^k z_l, y) )|^2}{|x-y|^{1+2s}} \,dx \,dy \\
 =
\int_{\mathbb{R}} \int_\mathbb{R}  \frac{| f( x') -  f(y')|^2}{|\chi(\sum_{l=1}^k z_l, x')) -\chi(\sum_{l=1}^k z_l, y') |^{1+2s}} \left|\frac{a(\chi(\sum_{l=1}^k z_l, x'))}{a(x')}\right| \left|\frac{a(\chi(\sum_{l=1}^k z_l, y'))}{a(y')}\right| \,dx' \,dy'. 
\end{multline*}
Now,  by the mean value theorem we have $|\chi(\sum_{l=1}^k z_l, x)) -\chi(\sum_{l=1}^k z_l, y)| \ge \frac{c_0}{c_1} |x-y|$ 
and the assumed bounds for $a$ give $\left|\frac{a(\chi(\sum_{l=1}^k z_l, \cdot))}{a(\cdot)}\right|\le \frac{c_1}{c_0}$, 
thus we arrive at
$$
\int_{\mathbb{R}} \int_\mathbb{R}  
  \frac{|g(x,z) - g(y,z)|^2}{|x-y|^{1+2s}} \,dx \,dy
 \leq
  \left(\frac{c_1}{c_0}\right)^{3+2s} 
  \int_{\mathbb{R}} \int_\mathbb{R}  \frac{| f( x) -  f(y)|^2}{|x-y|^{1+2s}} \,dx \,dy
$$
Again by $\int_{[0,\infty[^k}  |h(z)| \,d^k z = \frac{1}{{\rm Re}(\mu)^k}$  we conclude 
$$
 \|R(\mu)^k f \|_{s}^2 \le \frac{\frac{c_1}{c_0}}{{\rm Re(\mu)}^{2 k}} \cdot \| f\|_{0}^2   +
\frac{\left(\frac{c_1}{c_0}\right)^{3+2s}  }{{\rm Re(\mu)}^{2 k}} \int_{\mathbb{R}} \int_\mathbb{R}  \frac{| f( x') -  f(y')|^2}{|x'-y'|^{1+2s}} \,dx'  \,dy' 
\le
  \frac{\left(\frac{c_1}{c_0}\right)^{{3+2s}}  }{{\rm Re(\mu)}^{2 k}} 
    \|f\|_s^2,
$$
i.e., $\displaystyle{\|R(\mu)^k f \|_{s} \leq \frac{\left(\frac{c_1}{c_0}\right)^{(3+2s)/2}  }{{\rm Re(\mu)}^{k}}} \|f\|_s$.

\end{remark}



\subsection{Measurable coefficients with prescribed characteristics}

This subsection discusses  a solution concept according to Bouchut-James (\cite{BochutJames:98}), which is settled in one space dimension and --- from the distribution theoretic point of view --- can be considered as exotic.
The basic idea is to interpret the multiplication $a \cdot u$ occurring in the partial differential equation as a product of a (locally finite) Borel measure $u$ and a function $a$ from the set $\mathcal{B}^{\infty}$ of real bounded and Borel measurable functions.

\paragraph{Multiplication of Radon measures by bounded Borel functions:}
We may identify locally finite Borel measures on $\mathbb{R}$ with (positive) Radon-measures, that is the non-negative linear functionals on $C_c(\mathbb{R})$ (\cite[Remark 19.49]{HewStr:65}). Moreover, the space $\mathcal{D}^{'0}(\mathbb{R})$ is the space of complex Radon-measures, which allows for a decomposition of any
$u\in \mathcal{D}^{'0}(\mathbb{R})$ in the form $u= \nu_{+} -\nu_{-} + i (\eta_{+} - \eta_{-})$, where $\nu_{+},\nu_{-},\eta_{+},\eta_{+}$ are positive Radon-measures.

The product of a bounded Borel function $a\in  \mathcal{B}^{\infty}(\mathbb{R})$  with a positive Radon measure $\mu$ is defined to be the
measure given by
\begin{eqnarray*}
(a\odot \mu)(B) := \int_{\mathbb{R}^n} 1_B(x) a(x) \,d\mu(x) , \ \
\end{eqnarray*} 
for all Borel sets $B$ in $\mathbb{R}$. Clearly, $a\odot \mu$ is again a locally finite Borel measure.

The product employed in \cite{BochutJames:98} is the extension of $\odot$ to $ \mathcal{B}^{\infty}(\mathbb{R}) \times \mathcal{D}^{'0}(\mathbb{R})$ in a bilinear way, i.e.
\begin{eqnarray*}
\diamond: \mathcal{B}^{\infty}(\mathbb{R}) \times \mathcal{D}^{'0}(\mathbb{R}) &\rightarrow& \mathcal{D}^{'0}(\mathbb{R})\\
		(a, u) &\mapsto&   
a_{+}\odot \nu_{+} + a_{-}\odot \nu_{-} - ( a_{-}\odot \nu_{+}+a_{+}\odot \nu_{-}) \\
&&+i (a_{+}\odot\eta_{+} + a_{-}\odot \eta_{-} - ( a_{-}\odot \eta_{+}+a_{+}\odot \eta_{-}) ).
\end{eqnarray*}

Consider the following sequence of maps:
$$
C_b(\mathbb{R}) \stackrel{\iota_1}{\hookrightarrow} \mathcal{B}^{\infty}(\mathbb{R}) \stackrel{\lambda}{\rightarrow} L^{\infty}_{\rm loc}(\mathbb{R}) \stackrel{\iota_2}{\hookrightarrow} \mathcal{D}^{(0)\prime}(\mathbb{R})
$$
where $\iota_1, \iota_2$ are the standard embeddings and $\lambda$ sends bounded Borel functions to the corresponding classes modulo functions vanishing almost everywhere in
the sense of the Lebesgue measure.
%
%
Although we may identify $C_b(\mathbb{R})$ and $L^{\infty}_{\rm loc}$ with subspaces of $\mathcal{D}^{\prime 0}(\mathbb{R})$ this is not true of $\mathcal{B}^{\infty}(\mathbb{R})$, since $\lambda$ is not injective. Note that $\iota_2 \circ \lambda \circ \iota_1$ is injective though. The following example illustrates some consequences of the non-injectivity of the map $\lambda$ for the properties of the product $\diamond$.

\begin{example}
Let $\alpha\in\mathbb{R}$ and $a_{\alpha}(x):= 1,  x\ne 0$ and $a_{\alpha}(0):=\alpha$ and $u=\delta \in \mathcal{D}^{\prime0}(\mathbb{R})$.
 Note that $\lambda\circ \iota_1 (a_\alpha) =1$ as a distribution and the standard distributional product gives $\lambda\circ \iota_1(a_{\alpha}) \cdot u = 1 \cdot u = \delta$ for all $\alpha\in\mathbb{R}$. On the other hand $a_{\alpha}\diamond \delta = \alpha \delta$. 
\end{example}

The product $\diamond$ will be used in the solution concept for transport equations on $\overline{\Omega_T} = [0,T] \times \mathbb{R}$ with coefficient $a$
in $\mathcal{B}^{\infty}(\ovl{\Omega_T})$ and solution $u\in \mathcal{B}^{\infty}([0,T];\mathcal{D}^{'0}(\mathbb{R}))$, i.e., $u$ is a family of distributions $(u(t))_{t\in[0,T]}$ such
that $\langle u(t), \varphi\rangle_{\mathcal{D}^{'0}(\mathbb{R})}$ is a bounded Borel function on $[0,T]$ for all $\varphi \in C_c(\mathbb{R})$. The
extension of the product $\diamond$ to this space causes no difficulty.

\paragraph{The solution concept according to Bouchut and James:} A key ingredient for the solution concept according to Bouchut-James, is to stick to a particular representative of the coefficient (in the $L^{\infty}$ sense), by prescribing the value of the coefficient $a$ at curves of discontinuity. We refer to the following requirements on the coefficient $a\in \mathcal{B}^{\infty}(\Omega_T)$ as \emph{Bouchut-James conditions}: Assume there exists a decomposition $\Omega_T=\mathcal{C} \cup \mathcal{D} \cup \mathcal{S}$ such that
\begin{enumerate} \label{bochut_james_conditions}
\item $\mathcal{S}$ is a discrete subset $\Omega_T$,
\item $\mathcal{C}$ is open, $a$ is continuous on $\mathcal{C}$,
\item $\mathcal{D}$ is a one-dimensional $C^1$-submanifold  of ${\Omega_T}$,
i.e., for each $(t_0,x_0) \in \mathcal{D}$ there exists a neighborhood $V$ of $(t_0,x_0)$ and a $C^1$ parametrization of the form $t \mapsto (t,\xi(t))$ in $\D \cap V$.
Furthermore, $a$ has limit values for each $(t,x) \in \mathcal{D}$ from both sides in $\mathcal{C} \setminus \mathcal{D}$. These limits are denoted by $a_{+}(t,x)$ and $a_{-}(t,x)$.
\item $a(t,x) \in [a_{-}(t,x), a_{+}(t,x)]$ for all $(t,x) \in \mathcal{D}$,
\item for any point $(t_0,x_0) \in \mathcal{D}$ with neighborhood $V$ and local parametrization $\xi$ as in (iii), we have 
 $\dot{\xi}(t) = a(t,\xi(t)) $.
\end{enumerate}
Condition (v) prescribes the values of the coefficient $a(t,x)$ on the curves of discontinuity in such a way that
the characteristic differential equation holds. In this sense,  a coefficient satisfying $(i)-(v)$ is a piecewise continuous bounded function, where the (non-intersecting) curves of discontinuity can be parametrized as regular $C^1$ curves.

The Bouchut-James solution concept interprets hyperbolic Cauchy problems in $(1+1)$ dimension as
\begin{eqnarray} \label{pdo_P_bochut_james}
P u :=\partial_t u +  a \diamond \partial_{x}  u = 0, \ \ u(0)=u_0 \in {\rm BV}_{\rm loc}(\mathbb{R})  
\end{eqnarray}
and 
\begin{eqnarray} \label{pdo_L_bochut_james}
L u:= \partial_t u + \partial_{x} (a \diamond u)=0, \ \ u(0)=u_0 \in \mathcal{D}'^{ 0}(\mathbb{R}).
\end{eqnarray}

Note that $P$ (resp.\ L) is well-defined on the set $\mathcal{B}^{\infty}([0,T];{\rm BV}_{\rm loc}(\mathbb{R}))$ (resp.\ $\mathcal{B}^{\infty}([0,T]; \mathcal{D}^{\prime 0}(\mathbb{R}))$).

Now the main results of Bouchut-James \cite{BochutJames:98} are:

\begin{theorem}\cite[Theorem 3.4]{BochutJames:98}\label{BJexistence}
Assume that $a$ satisfies the Bouchut-James conditions (i)-(v).
For any $u_0\in BV_{\rm loc}(\mathbb{R})$ there exists $u\in {\rm Lip}([0,T]; L^1_{\rm loc}(\mathbb{R})) \cap \mathcal{B}^{\infty}([0,T]; BV_{\rm loc}(\mathbb{R}))$
solving (\ref{pdo_P_bochut_james}) and such that for any $x_1<x_2$ we have for all $t\in[0,T]$
\begin{eqnarray*}
{\rm Var}_I (u(t,\cdot)) &\le& {\rm Var}_J(u_0),\\
\|u(t,\cdot)\|_{L^{\infty}(I)} &\le& \|u_0\|_{L^{\infty}(J)}
\end{eqnarray*}
where $I:=]x_1,x_2[$ and $J:=]x_1-\|a\|_{\infty} t, x_2 + \|a\|_{\infty} t [$. If in addition the coefficient $a$ 
satisfies the one-sided Lipschitz condition
\begin{eqnarray*}
\langle a(t,x) -a(t,y), x-y \rangle \le \alpha(t) |x-y|^2 \ \ \text{for almost all} \ (t,x), (t,y) \in \Omega_T,
\end{eqnarray*}
where $\alpha \in L^1([0,T])$, then the solution $u$ is unique.
\end{theorem}

\begin{theorem}\cite[Theorem 3.6]{BochutJames:98} 
\label{BJexistence}
Assume that $a$ satisfies the Bouchut-James conditions (i)-(v).
Then it follows that for any $u_0\in \mathcal{D}^{'0}(\mathbb{R})$ there exists $u \in C([0,T]; \mathcal{D}^{'0}(\mathbb{R}))$
solving (\ref{pdo_L_bochut_james}). 
If $a$ satisfies in addition the one-sided Lipschitz condition
\begin{eqnarray*}
\langle a(t,x) -a(t,y), x-y \rangle \le \alpha(t) |x-y|^2 \ \ \text{for almost all} \ (t,x), (t,y) \in \Omega_T,
\end{eqnarray*}
where $\alpha \in L^1([0,T])$, then the solution $u$ is unique.
\end{theorem}

We compare the solution concept of Bouchut-James with the generalized solutions according to Poupaud-Rascle.
\begin{example} \label{BJ_jump} 
We come back to Example \ref{poupaud_rascle_jump}, where $a(t,x):=c_1 H(\alpha t- x) + c_2 H(x- \alpha t)$ with $c_2 < c_1$ and $\alpha\in [c_2,c_1]$.
Let $\lambda \in L^1([0,T])$ such that $c_1 \ge \lambda(t) \ge c_2$. 
Consider a representative $a$ in $\mathcal{B}^{\infty}(\overline{\Omega_T})$ of the coefficient $a$ given by
\begin{eqnarray*}
\widetilde{a}(t,x) =\left\{ 
\begin{array}{ll}
c_1 & x< \alpha t \\
\lambda(t) & x=\alpha t \\
c_2 &  x > \alpha t
\end{array}\right. .
\end{eqnarray*}

We investigate wether the distribution $u$ given in Example \ref{poupaud_rascle_jump} solves (\ref{pdo_L_bochut_james}) in the sense of Bouchut-James.

Let us consider the case $u_0\equiv 1$. Then we obtain the solution $u= 1+t (c_1-c_2) \delta(x-\alpha t)$.
Note that the requirement that $\widetilde{a}$ fulfills the {\emph Bouchut-James  conditions} (i)-(v) forces $\lambda(t) := \alpha$ for all $t\in [0,T]$. Thus we have 
$$
  \widetilde{a} \diamond v = 
   c_1  H(\alpha t-x)  + \alpha \delta(x-\alpha t)  t ( c_1 -c_2) + 
   c_2 H(x - \alpha t) 
$$
and therefore
$\partial_x (\widetilde{a} \diamond v) =   \alpha  t (c_1-c_2) \delta'(x-\alpha t) - (c_1-c_2) \delta(x-\alpha t)$. Since
$ \partial_t v =  (c_1-c_2) \delta(x-\alpha t) - \alpha t (c_1-c_2) \delta'(x-\alpha t)$ we deduce that
$u$ solves the Cauchy problem in the sense of Bouchut-James.

Finally, we check whether the differential equation is fulfilled, if we employ the model product (cf.\ \cite[Chapter 7]{O:92} or the introduction) instead of $\diamond$. 
Let $[a\cdot v]$ denote the model product. We have that
\begin{multline*}
  [\widetilde{a} \cdot v] = 
  [(c_1  H(\alpha t -x)+ c_2  H(x-\alpha t)) \cdot u] \\   
  = c_1 H(\alpha t -x) + c_2 H(x-\alpha t) + 
  c_1 t (c_1 -c_2) [H(\alpha t -x)   \cdot \delta(x-\alpha t)] \\
  = c_2 t (c_1 -c_2) [H(x-\alpha t) \cdot \delta(x-\alpha t)] 
+ c_1 H(\alpha t -x) + c_2 H(x-\alpha t) +  \frac{t}{2} (c_1+c_2) (c_1 -c_2)  \delta(x-\alpha t), 
\end{multline*}
hence
\begin{eqnarray*}
\partial_x [a \cdot v] &=& (c_2-c_1) \delta(x-\alpha t) - \frac{t}{2} (c_1+c_2) (c_1 -c_2)  \delta'(x-\alpha t),\\
\partial_t v &=&   (c_1-c_2) \delta(x-\alpha t) - \alpha t (c_1-c_2) \delta'(x-\alpha t).
\end{eqnarray*}
Therefore $v$ solves the initial value problem
\begin{eqnarray*}
\partial_t v + \partial_x [a\cdot v ] =0,\ \ v(0)=1,
\end{eqnarray*}
if $\alpha=\frac{1}{2}(c_1+c_2)$.

Note that the coefficients $H(-x)$ (with $c_1=1,c_2=0, \alpha=0$), $-H(x)$ (with $c_1=0,c_2=-1, \alpha=0$), and $-{\rm sign}(x)$ (when  $c_1=1,c_2=-1, \alpha=0$) are included as special cases of the example
presented here. In case the coefficient reads
$$
  a(x):= H(-x)
$$ 
the unique solution in the sense of Bouchut-James is given by 
$$ 
   u = 1+t \delta.
$$
It has been shown in \cite[Theorem 5]{HdH:01} that no distributional solution  exists in this case when the model product is employed.
\end{example}
%
%

%
%

\section{Solutions from energy estimates}

\subsection{Direct energy estimates}
We briefly review the standard techniques of energy estimates for the initial value problem 
\begin{equation} \label{initial}
\begin{split}
P u &:= \partial_t u +\sum_{j=1}^n a_j \,\partial_{x_j} u + c\, u = f \qquad \text{in } ]0,T[ \times \R^n, \\ 
u(0)&= u_0 \in L^2(\R^n). 
\end{split}
\end{equation}

Let $q \in [2,\infty]$. We assume that $f \in L^1([0,T];L^2(\R^n))$, 
$a = (a_1, ..., a_n) \in L^1([0,T]; W^{1,q}(\mathbb{R}^n))^n$ with real components,
$c \in  L^1([0,T]; L^q(\mathbb{R}^n))$
and in addition
\begin{equation} \label{div_cond} 
  \frac{1}{2} \, {\rm div}_x (a) -c \quad \in 
    L^1([0,T];L^{\infty}(\mathbb{R}^n)).
\end{equation}

\subsubsection{Example derivation of an energy estimate}

We browse through the typical steps that lead to an estimate in the norm of $L^{\infty}([0,T]; L^2(\mathbb{R}^n))$ for any 
\begin{eqnarray*}
u\in AC([0,T]; L^2(\mathbb{R}^n)) \cap L^{\infty}([0,T]; W^{1,p}(\mathbb{R}^n))
\end{eqnarray*}
with $p \in [2,\infty]$ such that $\frac{1}{q} + \frac{1}{p} = \frac{1}{2}$ in terms of corresponding norms for $u(0)$ and $Pu$.



We write $P = \d_t + Q$ with $Q:=\sum_{k=1}^n a_k(x,t) \partial_{x_k} + c(x,t)$ and observe that 
$$
  P u \in  L^{1}([0,T];L^2(\mathbb{R}^n))
$$ 
holds since $\partial_t u \in L^{1}([0,T];L^2(\mathbb{R}^n))$ and
$Qu \in L^{1}([0,T];L^2(\mathbb{R}^n))$ (the latter follows from the facts that $\d_{x_j} u(t,.) \in L^2$ and $L^q \cdot L^p \subseteq L^2$ when $2/p + 2/q = 1$).
Hence
$r \mapsto {\rm Re}(\langle (P u)(r) , u(r) \rangle_0)$ is defined and in $L^1([0,T])$. Furthermore, the map $t \mapsto \| u(t,\cdot)\|_0$ is continuous. 

We put 
$$
  h(r):= \|\frac{1}{2} \div_x (a(r,\cdot)) - 
    c(r,\cdot)\|_{\infty} \quad \text{ and } \quad
   \lambda(r) := 2 \int_0^r h(s) \,ds \geq 0 
   \qquad (r \in [0,T]).
$$ 
By assumption, $h \in L^1([0,T])$ and $\lambda \in AC([0,T])$. 

The standard integration by parts argument gives the G\r{a}rding-type inequality
\beq \label{Garding}
   \frac{1}{2} (\langle Qu(\tau),u(\tau) \rangle_0 + 
   \langle u(\tau), Qu(\tau) \rangle_0 ) =  
   {\rm Re}(\langle Qu(\tau), u(\tau) \rangle_0)  
   \ge - h(\tau) \| u(\tau) \|_0^2,
\eeq   
and thus
\begin{eqnarray*}
  \lefteqn{\int_0^\tau e^{-\lambda(r)} 
     {\rm Re}(\langle (P u)(r) , u(r) \rangle_0)\, dr} \\
   &=& \frac{1}{2} \int_0^\tau e^{-\lambda(r)} 
       \diff{r} \| u(r) \|_0^2\, dr + 
     \int_0^\tau e^{-\lambda(r)} 
     {\rm Re} \langle  (Q u)(r), u(r) \rangle_0 \, dr \\
   &\ge& \frac{1}{2} e^{-\lambda(\tau)} 
     \| u(\tau) \|_0^2 - \frac{1}{2} \| u(0) \|_0^2 - 
     \int_0^\tau  \underbrace{\left(h(r) - 
     \frac{\dot{\lambda}(r)}{2}\right)}_{= 0} 
     e^{-\lambda(r)} \| u(r) \|_0^2\, dr.
\end{eqnarray*}
Therefore
\begin{eqnarray*}
  e^{-\lambda(\tau)} \|  u(\tau) \|_0^2 &\le&  
  \| u(0) \|_0^2  + 
    2 \int_0^\tau e^{-\lambda(r)} \| (P u)(r) \|_0  
    \| u(r) \|_0 \,dr \\
  &\leq& \| u(0) \|_0^2 +  
    2 \!\!\sup_{r\in[0,\tau]}{ \Big( e^{-\lambda(r)/2} 
    \|   u(r) \|_0 \Big)}  \int_0^\tau e^{-\lambda(r)/2} 
    \| (P u)(r) \|_0 \,dr,
\end{eqnarray*}
where we may take the supremum over $\tau \in [0,t]$ on the left-hand side and thus replace $\tau$ by $t$ on the right-hand upper bound. A simple algebraic manipulation then gives
$$
  \left( \sup_{r\in[0,t]} \| e^{-\lambda(r)/2} 
     u(r)\|_0- \int_0^t e^{-\lambda(r)/2}  
     \|(P u)(r)\|_0   \,dr \right)^2  
  \le  \left(\| u(0) \|_0+    \int_0^t  
    e^{-\lambda(r)/2} \| (P u)(r) \|_0 \, dr\right)^2.
$$

Upon removing the squares and multiplying by $\exp(\lambda(t)/2)$ we obtain the following basic inequality.
\paragraph{Energy estimate:}
\begin{eqnarray}\label{energyest}
  \sup_{r\in[0,t]}\| u(r) \|_0 &\le&
   \exp (\int_0^t h(\sigma) d\sigma) \cdot 
   \| u(0) \|_0 +  2   \exp(\int_0^t h(\sigma) d\sigma)  
   \cdot \int_0^t  \| (P u)(r) \|_0 \,dr \\
   &=& \exp (\int_0^t h(\sigma) d\sigma) \left( 
   \| u(0) \|_0 + 2 \int_0^t  \| (P u)(r) \|_0 \,dr \right). \nonumber
\end{eqnarray}
We recall that the exponential factor depends explicitly on the coefficients $a$ and $c$ via $h(r) = 
\|\frac{1}{2} \div_x (a(r,\cdot)) - c(r,\cdot)\|_{\infty}$.

Note that this derivation of an energy estimate relied on the G\r{a}rding inequality \eqref{Garding}. 

\begin{example}[Failure of the G\r{a}rding-inequality \eqref{Garding}]
Let $\alpha\in \,]1/2,1[$ and define $a \col \R \to \R$ by $a(x) := 1+x_{+}^{\alpha}$ when $x\le 1$, and $a(x) := 2$ when $ x > 1$. 
 We have $a \in C_{\ast}^{\alpha}(\mathbb{R}) \setminus {\rm Lip}(\mathbb{R})$. 

Let $Q \col H^1(\R) \to L^2(\R)$ be the operator defined by $(Q v) (x):= a(x) v'(x)$ for all $v \in H^1(\R)$. Note that compared to the general form of the  operator $Q$ in the derivation of the energy estimate above we have here $c = 0$, $a \in C^{\infty}([0,T]; W^{1,2}(\R))$ but $\div a /2 - c = a'/2 \not\in L^{\infty}(\R)$.

Since $Q$ is time independent, inequality \eqref{Garding} with some $h \in L^1([0,T])$ (not necessarily of the form given above) would imply 
$$
  \exists C \in \R, \; 
  \forall v \in C^{\infty}_{\textrm{c}}(\R): \quad 
  {\rm Re} (\langle Q v, v \rangle_0)  
   \ge - C \| v \|_0^2.
$$ 
We will show that there is no constant $C \in \R$ such that the latter holds. Thus \eqref{Garding} cannot hold for $Q$ (for any $h \in L^1([0,T])$).

Let $\rho\in C^1(\mathbb{R})$ be symmetric, non-negative, with
support in $[-1,1]$, $\|\rho\|_0 =1$,  and such that $\rho'(x)<0$ when $0 < x < 1$.
We define $v_{\varepsilon}(x):= \varepsilon^{-1/2} \rho(x/\varepsilon)$  ($x \in \R$, $\eps > 0$). Then clearly $v_{\varepsilon} \in C^{\infty}_{\textrm{c}}(\R) \subseteq H^1(\mathbb{R})$ and $\|v_\varepsilon \|_0 = 1$ for all $\eps > 0$, but 
\begin{eqnarray*}
 \lefteqn{\langle Q v_{\varepsilon}, v_{\varepsilon} \rangle_0 =  
  \int a(x) v_{\varepsilon}'(x) v_{\varepsilon}(x)\, dx} \\
  &=& \underbrace{\int_{-\infty}^{\infty} v_{\varepsilon}'(x) 
  v_{\varepsilon}(x)\, dx}_{= 0} + \int_0^1 x^{\alpha} 
  v_{\varepsilon}'(x) v_{\varepsilon}(x)\, dx + 
   \underbrace{\int_1^{\infty} v_{\varepsilon}'(x) 
  v_{\varepsilon}(x)\, dx}_{=0} \\
  &=& \varepsilon^{\alpha - 1} \int_0^{1} 
  z^{\alpha} \rho'(z) \rho(z)\, dz \rightarrow - \infty 
  \qquad (\eps \to 0).
\end{eqnarray*}

We remark that even for $a \in C_{\ast}^{1}(\R) \setminus {\rm Lip}(\mathbb{R})$ the G\r{a}rding inequality may fail as well: for example, with $a(x):= - x \log{|x|} \rho(x)$ we have $a \in C_{\ast}^{1}(\R) \cap W^{1,q}(\R)$ for all $q \in [1,\infty[$, but 
$$
  \langle Q v_{\varepsilon}, v_{\varepsilon} \rangle_0 = 
   - 2 \int_0^1 \rho(\eps z) z \log |\varepsilon z|  
   \rho'(z) \rho(z) \, dz  
   \leq  2 |{\log{\varepsilon}}| \int_0^1 z  
   \rho(\eps z) \rho'(z) \rho(z) \, dz  \rightarrow - \infty,
$$
since $\lim_{\eps \to 0} \int_0^1 z \rho(\eps z) \rho'(z) \rho(z) \, dz = \rho(0) \int_0^1 z \rho'(z) \rho(z) \, dz < 0$.

\end{example}

\begin{remark}\label{adj_rem} (i) Let $Q^{\ast}$ denote the formal adjoint of $Q$ with respect to the $L^2$ inner product (on $x$-space). Due to our regularity assumptions on $a$ and $c$ we have for any $\vphi \in H^1$  (since $a$ is real)
$$
  Q^{\ast} \vphi = \sum_{j=1}^n (- a_j \d_{x_j} \vphi) 
  + (\bar{c} - \div_x(a)) \vphi,
$$
where the new coefficients $-a$, respectively $\bar{c} - \div_x(a)$, in place of $a$, respectively $c$, satisfy the exact same regularity assumptions, including the condition
$$
  \frac{1}{2} \div_x(- a) - (\bar{c} - \div_x(a)) =
  \ovl{\frac{\div_x(a) }{2} - c} \quad \in
   L^1([0,T];L^{\infty}(\mathbb{R}^n)).
$$
Thus the basic energy estimate \eqref{energyest} applies to $ \pm \d_t + Q^\ast$ as well. In particular, the function $h$ in the exponential factor occurring in the energy estimates is the same for $Q$ and $Q^\ast$.

(ii) Although the method of derivation discussed above relied on a G\r{a}rding-type inequality, it seems that in essence energy estimates are, in a vague sense, a necessary condition for a hyperbolic equation to hold in any meaningful context of ``suitable Banach spaces of distributions''. In other words, whenever a hyperbolic differential equation can be interpreted directly in terms of such Banach spaces it allows to draw consequences on combinations of corresponding norms of any solution. For example, if the operator $Q$ above generates a strongly continuous evolution system on some Banach space, then basic norm estimates for solutions follow from general principles of that theory (cf.\ \cite{Pazy:83,Tanabe:79}).
 \end{remark}

On the other hand, energy estimates are widely used to establish existence of solutions to \eqref{initial} by duality and an application of the Hahn-Banach theorem. We recall the basic steps of such method in the following.

\subsubsection{Existence proof based on the energy estimate}


Let $R_T:=\{(t,x)\in \mathbb{R}^{n+1} \mid  t < T\}$. By abuse of notation we denote the trivial extension of a function $v\in C_c^{\infty}(R_T)$ by zero for $t\ge T$
again by $v$. Then $\mathcal{L}:= \{ f \in C^{\infty}([0,T]\times \R^n) \mid  \exists v\in C_c^{\infty}(R_T) \text{ with }
f = (-\partial_t v+ Q^{\ast} v) \mid_{[0,T] \times \R^n} \}$. For $0  \le t \leq T$ and $v\in C_c^{\infty}(R_T)$ we use the notation 
$w(t):= v(T-t)$ and $g(t):= (-\partial_t v + Q^{\ast}v)(t)$.
Then we have
\begin{eqnarray*}
  (\partial_t + Q^{\ast}(T-t)) w(t) &=& g(T-t)\\
  w(0) &=&  0
\end{eqnarray*}
and an application of \eqref{energyest} (with $Q^*$ in place of $Q$; cf.\ Remark \ref{adj_rem}(i) above) yields 
\begin{multline*}
\sup_{r\in[0,T]} \| w(r) \|_{0} 
  \le  2 \exp{(\int_0^T h( \sigma ) \,d\sigma)} 
  \int_0^T  
  \| (-\partial_t+ Q^{\ast} v)(T-r) \|_{0}\, dr 
  = C_h \int_0^T \|g(r)\|_{0} \, dr.
\end{multline*}
We may deduce that for $f \in L^1([0,T]; L^2(\mathbb{R}^n))$ and $v\in C_c^{\infty}(R_T)$
\begin{multline*}
   \int_0^T \langle f(r), v(r) \rangle_0 dr + 
   \langle u_0, v(0) \rangle_0 \le 
   \int_0^T \|f(r)\|_{0} \|v(r)\|_{0} \,dr
   + \|u_0\|_{0}\, \|v(0)\|_{0}  \\
    \le C \sup_{r\in[0,T]}\|w(r)\|_{0} 
    \le C C_h \int_0^T \|g(r)\|_{0} \, dr,  
\end{multline*}
where $C$ depends on $f$ and $u_0$. Therefore the assignment
$g = (-\partial_t v + Q^{\ast} v) \mid_{[0,T] \times \R^n}\mapsto 
\int_0^T \langle f(r), v(r) \rangle_0 dr + \langle u_0, v(0) \rangle_0$ defines a conjugate-linear functional $\nu \col \mathcal{L} \to \mathbb{C}$ on the  subspace $\mathcal{L}$ of $L^1([0,T]; L^2(\mathbb{R}^n))$ such that
$|\nu(g)| \le \sup_{0 \leq r \leq T} \|g(r)\|_{0}$.  Hahn-Banach extension of $\nu$ yields a conjugate-linear functional
$\nu': L^1([0,T];L^2(\mathbb{R}^n)) \to \mathbb{C}$ with the same norm estimate. 

Since $L^1([0,T]; L^2(\mathbb{R}^n))' \cong L^{\infty}([0,T]; L^2(\mathbb{R}^n))$ there is $u \in L^{\infty}([0,T]; L^2(\mathbb{R}^n))$ such that 
$\nu'(g) = \langle u, g \rangle$ for all $g\in L^1([0,T]; L^2(\mathbb{R}^n))$. When applied to $g = (-\partial_t v + Q^{\ast} v) \mid_{[0,T] \times \R^n}$  with $v\in C_c^{\infty}(R_T)$ we obtain
\begin{multline} \label{DE_orig}
  \int_0^T \langle u(t),  -\partial_t v(t) + 
  ( Q^{\ast} v )(t) \rangle_0 \, dt =
  \dis{u}{((-\d_t + Q^{\ast}) v) \mid_{[0,T] \times \R^n}}  \\
  = \int_0^T \langle f(t), v(t,.)\rangle_0 \, dt + 
  \langle u_0, v(0) \rangle_0. 
\end{multline}

\subsubsection{Model discussion of the weak solution concept}

\paragraph{Case of smooth symbol:}
If the coefficients of $Q$ (and thus of $Q^{\ast}$) are $C^\infty$ then the above identity implies that $u$ is a distributional solution to the 
partial differential equation $P u = f$ in $]0,T[ \times \R^n$.
In fact, with $\vphi \in C_c^{\infty}(]0,T[ \times \R^n)$ in place of $v$ we have
$$
  \dis{(\d_t + Q)u}{\vphi} = \dis{u}{(-\d_t + Q^{\ast}) \vphi} 
  = \int_0^T \langle f(t), \phi(t,.)\rangle_0 \, dt = \dis{f}{\vphi}.
$$
Moreover, since $Q u \in L^{\infty}([0,T];H^{-1}(\R^n))$ the differential equation implies that
$$
  \d_t u = f - Q u \quad \in L^1([0,T];H^{-1}(\R^n))
$$
and thus $u \in AC([0,T];H^{-1}(\R^n))$. In particular, it makes sense to speak of the initial value $u(0) \in \D'(\R^n)$. Integrating by parts on the left-hand side of \eqref{DE_orig} (now reading \eqref{DE_orig} from right to left, and duality brackets in appropriate dual pairs of spaces) yields for any $v\in C_c^{\infty}(R_T)$
\begin{multline*}
   \int_0^T \langle f(t), v(t)\rangle_0 \, dt + 
  \langle u_0, v(0) \rangle_0 = 
  \int_0^T \underbrace{\dis{\d_t u (t) + Q u(t)}{v(t)}}_{\langle f(t), v(t,.)\rangle_0} \, dt - 
  \dis{u(T)}{\underbrace{v(T)}_{= 0}}_0 + \dis{u(0)}{v(0)}_0,
\end{multline*}
hence $u(0) = u_0$

Of course, uniqueness of the solution as well as more precise regularity properties can be deduced in case of $C^{\infty}$ coefficients: For any $s \in \R$, $f \in L^1([0,T]; H^s(\R^n))$, and $u_0 \in H^s(\R^n)$ the solution $u$ is unique in the space $C([0,T];H^s(\R^n))$ (cf.\ \cite[Theorem 23.1.2]{Hoermander:V3}).  

\paragraph{Case of non-smooth symbol:}
The weaker regularity assumptions made above imply $Q^{\ast} v \in L^1([0,T];L^2(\R^n))$ for all $v \in C_c^{\infty}(R_T)$. We may thus define $Q u \in \D'(]0,T[ \times \R^n)$ by putting 
$$
  \dis{Qu}{\vphi} := \inp{u}{Q^{\ast} \bar{\vphi}}_0 
     \qquad \forall \vphi \in C_c^{\infty}(]0,T[ \times \R^n).
$$ 
Then equation \eqref{DE_orig} can be read as an equation in $\D'(]0,T[ \times \R^n)$, namely
$$ 
    \dis{\d_t u + Q u}{\vphi} = \dis{f}{\vphi}  
    \qquad \forall \vphi \in C_c^{\infty}(]0,T[ \times \R^n).  
$$
Furthermore, we can again show that the inital datum is attained: Note that in $Q u = \sum a_j \d_{x_j} u$ each term can be interpreted as a multiplication of functions in $L^1([0,T];H^1(\R^n))$ with distributions in $L^{\infty}([0,T];H^{-1}(\R^n))$ (since $u \in L^{\infty}([0,T];L^2(\R^n))$) in the sense of the duality method (cf.\ \cite[Chapter II, Section 5]{O:92}). Applying Proposition 5.2 in \cite{O:92} to the spatial variables in the products then yields $Qu \in L^1([0,T]; W^{-1,1}(\mathbb{R}^n))$. Reasoning similarly as above, the differential equation then gives
$$
  \d_t u = f - Q u \quad \in L^1([0,T];W^{-1,1}(\R^n)),
$$
which implies $u \in AC([0,T];W^{-1,1}(\R^n))$ and further also that $u(0) = u_0$.    

Again higher regularity of $u$ with respect to the time variable, namely $u \in C([0,T];L^2(\R^n))$ can be shown by means of regularization and passage to the limit (e.g., similarly as in \cite[proof of Theorem 2.8]{BGS:07}).


\subsection{Regularization and energy estimates}

Several advanced theories make use of regularization techniques or concepts at crucial steps in their construction of solutions. Some of these theories succeed by regularization and a careful passage to the limit via energy estimates (as with Hurd-Sattinger and Di Perna-Lions theories presented below). Others even base their solution concept on a further generalization of the weak solution concept beyond distribution and measure spaces and still obtain existence of solutions essentially from asymptotic stability of energy estimates (cf.\ the Lafon-Oberguggenberger theory below). 
 
We introduce the following notation for partial differential operators that will be used in the sequel
\begin{eqnarray} 
  P u & := & \partial_t u + 
    \sum_{k=1}^n a_k \, \partial_{x_k}  u + c\, u 
    \label{pdo1} \\
  L u & := & \partial_t u + 
   \sum_{j=1}^n \partial_{x_j} (a_j\, u) + b\, u.
\end{eqnarray}

\subsubsection{Hurd-Sattinger theory}

We give a brief summary of the results from the first part in Hurd-Sattinger's classic paper \cite{HS:68}. We consider the Cauchy problem for the operator $L$ on the closure of the domain $\Omega := \, ]0,\infty[ \times \R^n$.

\begin{definition} Let $f \in L^2_\loc(\ovl{\Omega})$ and $a_j$ ($j=1,\ldots,n$) as well as $b \in L^2_\loc(\ovl{\Omega})$. A weak solution in the sense of Hurd-Sattinger of the partial differential equation
$$
   L u = f \qquad \text{on } \Omega
$$
with initial condition $u_0 \in L^2_{\rm{loc}}$ is a function  
$u \in L^2_{\rm{loc}}(\ovl{\Omega})$ such that for all $\phi \in C_c^{1}(\mathbb{R}^{n+1})$ we have
\begin{multline}
   \int_{\ovl{\Omega}} \Big(- u(t,x) \phi(t,x) - 
  \sum_{j=1}^n 
  a_j(t,x) u(t,x) \partial_{x_j} \phi(t,x) 
  + b(t,x) u(t,x) \phi(t,x) \Big) \, d(t,x) \\
  = \int_{\ovl{\Omega}} f(t,x) \phi(t,x) \, d(t,x)
 + \int_{\mathbb{R}^n} u_0(x) \phi(0,x) \,dx.
\end{multline}
\end{definition}

Note that if all coefficients are  $C^{\infty}$ functions then a solution in the above sense solves the partial differential equation on $\Omega$ in the sense of distributions.

\begin{theorem}\label{HS_thm} Let $a_j$ ($j=1,\ldots,n$), $b$, and $f$ belong to $L^2_\loc(\ovl{\Omega})$ and $u_0 \in L^2_{\rm{loc}}$.
Assume, in addition, that the following conditions are satisfied:
\begin{enumerate}

\item There exists $c_1 > 0$ such that for almost all $(t,x) \in \ovl{\Omega}$: \enspace $a_k (t,x) \le c_1$ \enspace 
($k=1\ldots,n$).

\item There exists a function $\mu \in L^1_{\rm{loc}}([0,\infty[)$, $\mu \geq 0$, such that $b(t,x) \ge - \mu(t)$ for almost all $(t,x) \in \ovl{\Omega}$.

\item For each $k \in \{1,\ldots,n\}$ there exists $0\le\mu_k \in L^1_{\rm{loc}}([0,\infty[)$ such that for almost all $(t,x) \in \ovl{\Omega}$
$$ 
\frac{a_k(t,x) - a_k(t,x_1,\ldots,x_{k-1},r,x_{k+1}\ldots,x_n)}{x_k - r} \ge - \mu_k(t) \qquad \text{for almost all } r \in \R.
$$
\end{enumerate}
Then there exists a weak solution $u \in L^2_{\rm{loc}}(\ovl{\Omega})$ to $L u = f$ with initial condition $u_0$.
\end{theorem}

Concerning the meaning of condition (iii) in Theorem \ref{HS_thm} we mention two aspects:
\begin{itemize}
\item In one space dimension we obtain
$\frac{a(x,t)-a(t,y)}{x-y} \ge - \mu_1(t)$,
which resembles a one-sided Lipschitz continuity condition in the $x$ variable (apart from the fact that $\mu_1(t)$ need not be finite or defined for all $t$). In particular, it excludes jumps downward (seen when going from smaller to larger values in the $x$ argument).

\item Heuristically --- replacing difference quotients by partial derivatives --- condition (iii) can be read as $\div a(t,x) \geq - \sum \mu_k(t)$, thus giving an $L^1$ lower bound on the divergence of $a$. We observe that upon formally applying the Leibniz rule in the operator $L$ we cast it in the form $P$ as in \eqref{pdo1} with  $c = \div a + b$. In combination with condition (ii) of Theorem \ref{HS_thm}, we obtain that $\frac{1}{2} \div a - c =  -(\frac{1}{2} \div a + b)$ has an $L^1$ upper bound (uniformly in $x$), which can be considered a substitute for condition \eqref{div_cond} used in the derivation of direct energy estimates in Subsection 2.1.
\end{itemize}  

\begin{remark} Hurd-Sattinger (\cite{HS:68}) also give a uniqueness result for first-order systems in case of a single space variable and $b = 0$. For scalar equations the hypotheses require condition (i) to be strengthened to boundedness from above and from below and  condition (iii) to be replaced by a Lipschitz property with an upper bound instead; in particular, no jumps upward are possible.  
\end{remark}

\begin{example} For the operator $L$ in one space dimension and coefficients $a(x) = \rm{sign}(x)$ and $b = 0$, the Poupaud-Rascle theory is not applicable (as mentioned in \cite[Section 1, Example 2]{PoupaudRascle:97}), but Hurd-Sattinger theory ensures existence of weak solutions, if the initial value belongs to $L^2_{\rm{loc}}$.
\end{example}


\subsubsection{Di Perna-Lions theory}

The weak solution concept introduced by Di Perna-Lions in \cite{DiPernaLions:89} for the Cauchy problem for the operator $P$ on a finite-time domain $[0,T] \times \R^n$ can be interpreted in the following way.

\begin{definition}\label{dPL_sol} Let $T > 0$, $1 \leq p \leq \infty$, $\frac{1}{p} + \frac{1}{q} = 1$, $f \in L^1([0,T]; L^p(\R^n))$, $a_k \in L^{1}([0,T];L^q_{\rm{loc}}(\mathbb{R}^n))$ ($k=1,\ldots,n$), and $c\in L^1([0,T];L^q_{\rm loc}(\mathbb{R}^n))$ such that 
$$
  \div(a) - c  \quad \in 
  L^{1}([0,T];L^q_{\rm{loc}}(\mathbb{R}^n)).
$$ 
A function $u\in L^{\infty}([0,T]; L^{p}(\mathbb{R}^n))$ is called a weak solution in the sense of Di Perna-Lions  of the partial differential equation 
$$
 P u = f \qquad \text{on } ]0,T[\, \times \R^n
$$
with initial value $u_0 \in L^p(\R^n)$, if 
\begin{multline} 
  \int_0^T \int_{\mathbb{R}^n} u(t,x) \Big( -
  \partial_t \varphi(t,x) dx  
  - \sum_{k=1}^n a_k(t,x) 
  \d_{x_k}\varphi(t,x) \Big) \,dx dt \\ 
  + \int_0^T \int_{\mathbb{R}^n}
  u(t,x) \big( -
   \div a(t,x) + c(t,x) \big) 
  \varphi(t,x)  \,dx dt \\
  = \int_0^T \int_{\mathbb{R}^n} f(t,x) 
  \varphi(t,x)\, dx dt  + 
  \int_{\mathbb{R}^n} u_0(x) \varphi(0,x) \, dx 
\end{multline}
holds for all $\varphi \in C^{\infty}([0,T], \mathbb{R}^n)$ with compact support in $[0,T[ \times \R^n$.
\end{definition}

Clearly, in case of $C^{\infty}$ coefficients we obtain a distributional solution of the partial differential equation in $]0,T[ \times \R^n$.

\begin{theorem}\label{DiPL_thm} Existence of a weak solution $u \in L^{\infty}([0,T]; L^p(\mathbb{R}^n))$ in the sense of and with assumptions as in Definition \ref{dPL_sol} is guaranteed under the additional hypothesis 
\begin{eqnarray*}
  \frac{1}{p} {\rm div}(a) - c \in \quad
     L^1([0,T];L^{\infty}(\R^n)),
  && \rm{if}\ p>1, \\ 
    {\rm div}(a), c \quad \in 
    L^1([0,T];L^{\infty}(\R^n)),
    && \rm{if}\ p=1.
\end{eqnarray*}
\end{theorem}

\begin{remark}\label{DiPL_uni} Uniqueness holds in general under the
 additional hypotheses that $c, {\rm div}(a) \in L^1([0,T]; L^{\infty}(\mathbb{R}^n))$, and for  $j=1,\ldots,n$ also $a_j \in L^1([0,T]; W^{1,q}_{\rm loc}(\mathbb{R}^n))$ as well as
\begin{eqnarray*}
  \frac{a_j}{1+|x|} \quad \in \;
  L^1([0,T]; L^1(\mathbb{R}^n)) + 
  L^1([0,T]; L^{\infty}(\mathbb{R}^n)).
\end{eqnarray*}
\end{remark}

\begin{example}[Hurd-Sattinger applicable, but not Di Perna-Lions] Note that with a single spatial variable boundedness of $\div(a) = a'$ implies Lipschitz continuity. Hence, if $a \in H^1(\R)$ is not Lipschitz continuous but satisfies the one-sided Lipschitz condition in Hurd-Sattinger's existence Theorem \ref{HS_thm} (iii), then a weak solution in the sense of Hurd-Sattinger to the problem
$$
   \d_t u + \d_x(a u) = f \in L^2(\R^2), \quad
   u \mid_{t=0} = u_0 \in L^2(\R)
$$   
is guaranteed to exist, whereas the general statement of DiPerna-Lions' existence theory (Theorem \ref{DiPL_thm} with $p=q=2$) is not applicable to the formally equivalent problem
$$
  \d_t u + a \d_x u + a' u = f \in L^2(\R^2), \quad
   u \mid_{t=0} = u_0 \in L^2(\R).
$$ 
\end{example}

\begin{example}[Di Perna-Lions applicable, but not Hurd-Sattinger] 
Let $0 < \sigma < 1$ and consider the identical coefficient functions $a_1 = a_2 \in C_{*,\mathrm{comp}}^{\sigma}(\R^2)$ (i.e., compactly suported functions in $C_*^\sigma(\R^2)$) given by
$$
  a_1(x,y) = a_2(x,y) = 
  -\frac{1}{\sigma} (x-y)_{+}^{\sigma}\, \chi(x,y),
$$
where $\chi \in \mathcal{D}(\mathbb{R}^2)$ such that $\chi = 1$ near $(0,0)$. Note that $a_1$ is not Lipschitz continuous, since for $x > 0$ but $x$ sufficiently small the difference quotient
$$
  \frac{a_1(x,0) - a_1(0,0)}{x} = - \frac{x^{\sigma - 1}}{\sigma}
$$ 
is unbounded as $x \to 0$. In particular, the latter observation shows that the Hurd-Sattinger existence theory is not applicable (condition (iii) in Theorem \ref{HS_thm} is violated) to the Cauchy problem for the operator
$$
   L u = \d_t u + \d_x(a_1 u) + \d_y(a_2 u).
$$

On the other hand, we can show that with $a = (a_1,a_2)$ the DiPerna-Lions existence theory is applicable to the Cauchy problem
\begin{eqnarray*}
  \partial_t u + a_1 \d_x u + a_2 \d_y u  + (\div a)\, u 
  = f \in L^1([0,T]; L^p(\mathbb{R}^2)), 
  \quad u \mid_{t=0} = u_0 \in L^p(\R^2).
\end{eqnarray*} 

To begin with, we observe that
\begin{align*}
  \d_x a_1(x,y) = \d_x a_2(x,y) &=
  - \frac{\chi(x,y)}{(x-y)_+^{1-\sigma}} - 
  \frac{1}{\sigma} (x-y)_+^\sigma \, \d_x \chi(x,y) \\
  \d_y a_1(x,y) = \d_y a_2(x,y) &=
  \frac{\chi(x,y)}{(x-y)_+^{1-\sigma}} - 
  \frac{1}{\sigma} (x-y)_+^\sigma \, \d_y \chi(x,y)  
\end{align*} 
yields
$$
 \div a (x,y) = - 
  \frac{1}{\sigma} (x-y)_+^\sigma \, \div \chi(x,y) 
  \in C_{*,\mathrm{comp}}^\sigma(\R^2).
$$
Note that in the notation of Definition \ref{dPL_sol} and Theorem \ref{DiPL_thm} we have $c = \div a \in L^\infty(\R^2)$ (and time-independent). Therefore, the basic assumptions for the solution concept to make sense as well as the hypotheses of the existence statement are clearly satisfied.

As for uniqueness, we remark that all the conditions mentioned in Remark \ref{DiPL_uni} are met if and only if $\sigma > 1/p$.
\end{example}

\begin{remark} We mention that with coefficients as in the above example, the system of characteristic differential equations has  forward-unique solutions, hence the Poupaud-Rascle solution concept for measures is also applicable. 
\end{remark}


\subsubsection{Lafon-Oberguggenberger theory}

The theory for symmetric hyperbolic systems presented in \cite{LO:91} by Lafon-Oberguggenberger allows for Colombeau generalized functions as coefficients as well as inital data and right-hand side. Thus we consider the following hyperbolic Cauchy
problem in $\R^{n+1}$
\begin{align}
  P u = \d_t u + \sum_{j=1}^{n} a_j \d_{x_j} u + c u = f   
\label{hypsys_equ} \\
  u \mid_{t = 0} = u_0,
\label{hypsys_ini}
\end{align}
where $a_j$ ($j=1,\ldots,n$), $c$ are real valued generalized
functions in $\G(\R^{n+1})$ (in the sense that all representatives are
real valued smooth functions), $f \in \G(\R^{n+1})$, and initial value $u_0 \in\G(\R^n)$. 

The coefficients will be subject to some restriction on the allowed
divergence in terms of $\eps$-dependence. A Colombeau function
$v \in\G(\R^d)$ is said to be of \emph{logarithmic type} if it has a
representative $(v_\eps)$ with the following property: there are constants
$N\in\N$, $C>0$, and $1 > \eta>0$ such that
\begin{equation*}
  \sup\limits_{y\in\R^d} |v_\eps(y)| \leq N \log\big(\frac{C}{\eps}\big) 
  \qquad 0 <  \eps < \eta \; .
\end{equation*}
(This property then holds for any representative.) By a suitable modification of \cite[Proposition 1.5]{O:89}
it is always possible to model any finite order distribution as coefficient with such properties (in the sense that the Colombeau coefficient is associated to the original distributional coefficient).

\begin{theorem}\label{LO_thm} Assume that $a_j$ and $c$ are constant for large $|x|$ and that 
$\d_{x_k} a_j$ ($k=1,\ldots,n$) as well as $c$ are of logarithmic type. Then given initial data $u_0 \in\G(\R^n)$ and right-hand side $f \in \G(\R^{n+1})$, the Cauchy problem (\ref{hypsys_equ})-(\ref{hypsys_ini}) has a unique solution $u \in\G(\R^{n+1})$.
\end{theorem}

We also mention the following {\bf consistency result} which shows that
Colombeau theory includes the classically solvable cases: If we assume that the coefficients $a_j$ and $c$ are $C^\infty$ then we have the following consistency with classical and distributional solutions (cf.\ \cite{LO:91})
\begin{itemize}
\item If $f$ and $u_0$ are $C^\infty$ functions then the generalized solution
  $u \in\G(\R^{n+1})$ is equal (in $\G$) to the classical smooth
  solution.
\item If $f \in L^2(\R;H^s(\R^n))$ and $u_0 \in H^s(\R^n)$ for some
  $s\in\R$, then the generalized solution $u \in\G(\R^{n+1})$ is
  associated to the classical solution belonging to $C(\R;
  H^s(\R^n))$.
\end{itemize}

\begin{example} Consider the $(1+1)$-dimensional operator
$$
   L u  = \d_t u  + \d_x (H(-x) u).
$$
Since the coefficient (of the formal principal part) has a jump downward neither Hurd-Sattinger nor Di Perna-Lions theory is applicable. In fact, it has been shown in \cite[Section 2]{HdH:01} that none of the distributional products from the coherent hierarchy (cf.\ \cite{O:92} and the introductary section) applied to $H(-x) \cdot u$ is capable of allowing for distributional solutions of the homogeneous Cauchy problem for arbitrary smooth initial data.

Recall from Section 1 that measure solutions according to Bouchut-James exist for the corresponding Cauchy problem, if the Heaviside function (usually understood as a class of functions in $L^\infty$) is replaced by the particular Borel measurable representative with value $0$ at $x = 0$.  For example, the initial value $u_0 = 1$ then yields the measure solution $u = 1 + t \delta(x)$ in the sense of Bouchut-James as seen in Example \ref{BJ_jump}. 

However, Colombeau generalized solutions are easily obtained --- even for arbitrary generalized initial data --- if the coefficient $H(-x)$ is regularized by convolution with a delta net of the form $\rho_\eps(x) =  \log(1/\eps) \rho(x \log(1/\eps))$ ($0 < \eps < 1$), where $\rho \in C_c^\infty(\R)$ with $\int \rho = 1$. Let $a$ denote the class of this regularization in the Colombeau algebra $\G$, then the operator $L$ may now be written equivalently in the form
$$
  P u = \d_t u + a \d_x u + a' u,
$$
where $a' \approx \delta$ and $u \in \G$. Due to the logarithmic scale in the regularization the hypotheses of Theorem  \ref{LO_thm} are satisfied and the corresponding Cauchy problem is uniquely solvable. Moreover, for most interesting initial data (e.g.\  Dirac measures or $L^1_{\mathrm{loc}}$) weak limits of the Colombeau solution $u$ are known to exist and can be computed (cf.\ \cite[Section 6]{HdH:01}). In particular, for the initial value $u_0 = 1$ we obtain the measure solution $u = 1 + t \delta(x)$ as such a distributional shadow.

\end{example}

\begin{remark} (i) The basic results of Lafon-Oberguggenberger have been extended to the case of (scalar) pseudodifferential equations with generalized symbols in \cite{GH:04}.  Special cases and very instructive examples can be found in \cite{O:88}, and an application of Colombeau theory to the linear acoustics system is presented in \cite{O:89}).

(ii) Colombeau-theoretic approaches allow for a further flexibility even in interpreting distributional differential equations with smooth coefficients. For example, in \cite{CHO:96} the concept of regularized derivatives is used, where partial differentiation is replaced by convolution with the corresponding derivative of a delta sequence. When acting on distributions this concept produces the usual differential operator actions in the limit. When considered as operators in Colombeau spaces, one can prove (cf.\ \cite[Theorem 4.1]{CHO:96}) that evolution equations with smooth coefficients all whose derivatives are bounded have unique generalized function solutions for initial data and right-hand side in generalized functions.
In particular, famous examples like the Lewy equation become solvable and Zuily's non-uniqueness examples become uniquely solvable then.
\end{remark}

\subsection{Paradifferential techniques}

\subsubsection{Energy estimates }

Bony's paradifferential calculus has been successfully applied  in nonlinear analysis and, in particular, to regularity theory for  nonlinear partial differential equations. An ingredient in such approaches is often a refined regularity assessment of corresponding linearizations of the differential operators involved. A recent account of M\'{e}tivier's methods and results of this type can be found in \cite[Subsection 2.1.3]{BGS:07}, or with more details on microlocal properties in \cite{Hoermander:97}.

Let $s \in \R$ and $H^s_w(\R^n)$ denote the Sobolev space $H^s(\R^n)$ equipped with the weak topology. We consider a differential operator of the form
$$
  \widetilde{P}_v (x,t;\d_t,\d_x) := 
  \d_t + \sum_{j=1}^n a_j(v(x,t)) \, \d_j,
$$
where $a_j \in C^\infty(\R)$ ($j=1\ldots,n$) and $v \in L^\infty([0,T];H^s(\R^n)) \cap C([0,T];H^s_w(\R^n))$ such that $\d_t v \in L^\infty([0,T];H^{s-1}(\R^n)) \cap C([0,T];H^{s-1}_w(\R^n))$. 

\begin{remark}
Not all hyperbolic first-order differential operators with coefficients of regularity as above can be written in the special form of $\widetilde{P}_v$. In fact, this amounts to writing any given list $w_1, \ldots, w_n$ of such functions as $w_j = a_j \circ v$ ($j=1\ldots,n$) with $a_j \in C^\infty(\R)$ and $v$ as above. The latter is, in general, not possible, which can be seen from the following example: consider the Lipschitz continuous functions $w_1(t) = |t|$ and $w_2(t) = t$; if $w_1 = a_1 \circ v$ and $w_2 = a_2 \circ v$ with a Lipschitz continuous function $v$, then $v$ is necessarily non-differentiable at $0$; on the other hand 
$$
  1 = w_2'(0) = \lim_{h \to 0} (a_2(v(h)) - a_2(v(0)))/h 
   = \lim_{h \to 0} a_2'(\xi(h)) (v(h) - v(0)) /h,
$$ 
where $\xi(h)$ lies between $v(0)$ and $v(h)$; hence $a_2'(\xi(h)) \to a_2'(v(0))$ and the second factor $(v(h) - v(0)) /h$ stays bounded, but is not convergent; in case $a_2'(v(0)) = 0$ we obtain the contradiction $1 = 0$, in case $a_2'(v(0)) \neq 0$ we have a contradiction to convergence of the difference quotient for $w_2$.   
\end{remark}

The key technique in analyzing the operator $\widetilde{P}_v$ is to replace all terms $a_j(v)  \d_j$ by $T_{a_j(v)}  \d_j$, i.e., partial differentiation followed by the para-product operator $T_{a_j(v)}$, and then employ estimates of the error terms as well as a paradifferential variant of G\r{a}rding's inequality  (cf.\ \cite[Appendix C.3-4]{BGS:07}). This leads to the following result.
\begin{theorem}[{\cite[Theorem 2.7]{BGS:07}}]
If $s > \frac{n}{2} + 1$, then for any $f \in L^\infty([0,T];H^s(\R^n)) \cap C([0,T];H^s_w(\R^n))$ and $u_0 \in H^s(\R^n)$ the Cauchy problem
$$
  \widetilde{P}_v u = f, \quad u \mid_{t=0} = u_0 
$$
has a unique solution $u \in L^2([0,T];H^s(\R^n))$. Moreover, $u$ belongs to $C([0,T];H^s(\R^n))$ and there are constants $K, \gamma, C \geq 0$ such that $u$ satisfies the energy estimate
$$
  \norm{u(t)}{s}^2 \leq K e^{\gamma t} \norm{u(0)}{s}^2 
    + C \int_0^t e^{\gamma (t - \tau)} 
    \norm{\widetilde{P}_v u(\tau)}{s}^2 \, d\tau.
$$
\end{theorem}

\subsubsection{Improvement of regularity in one-way wave equations}

We briefly recall some basic notions and properties concerning symbols with certain H\"older regularity in $x$ and smoothness in $\xi$  \`{a} la Taylor (cf.\ \cite{Taylor:91}). 

\begin{definition} Let $r > 0$, $0 < \delta < 1$, and $m \in \R$. A continuous function $p: \R^n \times \R^n \to \C$ belongs to the symbol space $C_\ast^r S^m_{1,\delta}$, if for every fixed $x \in \R^n$ the map $\xi \mapsto p(x,\xi)$ is smooth and 
for all $\al \in \N_0^n$ there exists $C_\al > 0$ such that
$$
  |\d_\xi^\al p(x,\xi)| \leq C_\al (1 + |\xi|)^{m - |\al|}
  \qquad \forall x, \xi \in \R^n
$$
and
$$
  \norm{\d_\xi^\al p(.,\xi)}{C_\ast^r} \leq 
    C_\al (1 + |\xi|)^{m - |\al| + r \delta}
    \qquad \forall \xi \in \R^n.
$$
\end{definition}

Basic examples are, of course, provided by symbols of differential operators $\sum a_\al \d^\al$ with coefficient functions $a_\al \in C_\ast^r$ ($|\al| \leq m$) or any symbol of the form $p(x,\xi) = a(x) h(x,\xi)$, where $a \in C_\ast^r$ and $h$ is a smooth symbol of order $m$.

\paragraph{Symbol smoothing:} By a coupling of a Littlewood-Paley decomposition in $\xi$-space with convolution regularization in $x$-space via a $\delta$-dependent scale one obtains a decomposition of any symbol $p \in C_\ast^r S^m_{1,\delta}$ in the form 
$$
 p = p^\sharp + p^\flat, \quad \text{ where } p^{\sharp} 
     \in S^m_{1,\delta} \text{ and } p^{\flat} 
     \in C_\ast^r S^{m - r \delta}_{1,\delta}.
$$
Observe that $p^\sharp$ is $C^\infty$ and of the same order whereas $p^\flat$ has the same regularity as $p$ but is of lower order.

\paragraph{Mapping properties:}
Let $0 < \delta < 1$ and $-(1-\delta) r < s < r$. Then any symbol  
$p\in C_\ast^r S^{m}_{1,\delta}$ defines a continuous linear operator
$p(x,D): H^{s+m}(\mathbb{R}^n) \rightarrow H^s(\mathbb{R}^n)$.

\paragraph{Elliptic symbols:} $p \in C_\ast^r S^m_{1,\delta}$ is said to be elliptic, if there are constants $C, R > 0$ such that
$$
   |p(x,\xi)| \geq C (1 + |\xi|)^m 
     \qquad \forall \xi \in \R^n, |\xi| \geq R.
$$

One-way wave equations result typically from second-order partial differential equations by a pseudodifferential decoupling into two first-order equations (cf.\ \cite[Section
IX.1]{Taylor:81}). For example, this has become a standard technique in mathematical geophysics for the decoupling of modes in seismic wave propagation (cf.\ \cite{SdH:02}). The corresponding Cauchy problem with seismic source term $f \in C^{\infty}([0,T]; H^{s}(\mathbb{R}^n))$ (with $s \in \R$) and initial value of the displacement $u_0 \in H^{s+1}(\mathbb{R}^n)$ is of the form
\begin{eqnarray}
  \partial_t u + i \, Q(x,D) u &=& f \label{one-way-PDE}\\
  u\mid_{t=0} &=& u_0 \label{one-way-initial},
\end{eqnarray}
where $Q$ has real-valued elliptic symbol $q \in C_\ast^r S^1$ with $r > s$.

\begin{lemma}\label{ell_lemma}
If $q\in C^{r}S^m_{1,0}$ is elliptic, then $q^{\sharp} \in S^{m}_{1,\delta}$ is also elliptic. 
\end{lemma}
\begin{proof}
By ellipticity of $q$ and the symbol properties of $q^\flat$ there are constants $C_1, C_2, R > 0$ such that 
\begin{eqnarray*}
C_1(1+|\xi|)^m \le |q(x,\xi)|\le |q^{\sharp}(x,\xi)| +|q^{\flat}(x,\xi)| \le  |q^{\sharp}(x,\xi)| + C_2 (1+|\xi|)^{m-r\delta}
\end{eqnarray*}
holds for all $x, \xi \in \R^n$ with $|\xi| \geq R > 0$. 
Therefore
$$
|q^{\sharp}(x,\xi)| \ge 
  (C_1 - C_2 (1 + |\xi|)^{-r\delta}) (1+|\xi|)^m 
  \geq C (1 + |\xi|)^m  
  \qquad \forall x, \xi \in \R^n, |\xi| \geq R'
$$
for suitably chosen constants $C$ and $R' > 0$.
\end{proof}

Let $0 < \delta < 1$. We have the decomposition $q = q^{\sharp} + q^{\flat}$, where
$q^{\sharp} \in S^1_{1,\delta}$ and $q^{\flat} \in C^{r} S^{1-\delta r}_{1,\delta}$. By Lemma \ref{ell_lemma} $Q^{\sharp} = q^\sharp(x,D)$ is elliptic and thus possesses a parametrix $E^{\sharp} \in S^{-1}_{1,\delta}$.

We have\
$$
  (\partial_t + i Q) E^{\sharp} f = 
  (\partial_t + i Q^\sharp + i Q^\flat) E^{\sharp} f 
  = \partial_t E^{\sharp} f + i Q^{\sharp} E^{\sharp} f 
    + i Q^{\flat} E^{\sharp} f
  =  \partial_t E^{\sharp} f + f +i R^{\sharp} f +
    i Q^{\flat} E^{\sharp} f,
$$
where $R^{\sharp}$ is a regularizing operator.
Therefore
$$
   (\partial_t + i Q) (u - E^{\sharp} f) = 
   -\partial_t E^{\sharp} f  -i R^{\sharp} f - 
   i Q^{\flat} E^{\sharp} f =: \widetilde{f},
$$
where the regularity of the right-hand side $\widetilde{f}$ can be deduced from the following facts
$$
  \partial_t E^{\sharp} f \in 
    C^{\infty}([0,T]; H^{s+1}(\mathbb{R}^n)), \quad
 R^{\sharp} f \in C^{\infty}([0,T]; H^{\infty}(\mathbb{R}^n)),
 \quad
 Q^{\flat} E^{\sharp} f \in C^{\infty}([0,T]; H^{s +\delta r}(\mathbb{R}^n)).
$$
Hence $\widetilde{f} \in C^{\infty}([0,T]; H^{s + \min{(\delta r,1)}}(\mathbb{R}^n))$.

If we put  
$w= u+ E^{\sharp} f$  and $w_0 := u_0 + E^{\sharp} f(0)$, 
then the original Cauchy problem (\ref{one-way-PDE}-\ref{one-way-initial}) is reduced to solving the Cauchy problem
$$
  \partial_t w + i\, Q(x,D) w = \widetilde{f}, 
  \qquad w \mid_{t=0} = w_0,
$$
where the spatial regularity of the source term on the right-hand side has been raised by $\min (\delta r,1)$.

\begin{remark} In case of a homogeneous ($1+1$)-dimensional partial differential equation the precise H\"older-regularity properties  
of classical as well as generalized solutions
have been determined in \cite[Section 3]{GH:04b}.
\end{remark}

\paragraph{Acknowledgements:} We are deeply thankful to Michael Oberguggenberger for having introduced the second author to several of the topics discussed in this paper and for numerous joint discussions in Vienna and Innsbruck. In particular, we are indebted to him for sharing with us his initial observations on semigroups generated by operators with H\"older-continuous coefficients, which gave rise to the examples presented in Subsection 1.4.

\bibliographystyle{abbrv}
\bibliography{simhal}

\end{document}